\documentclass[11pt]{amsart}

\usepackage{a4}
\usepackage{amsmath}
\usepackage{amsfonts}
\usepackage{amsthm}
\usepackage{amssymb}
\usepackage{times}
\usepackage{color}
\usepackage{graphicx}
\usepackage{url}

\newcommand{\finishline}{\hfill\hbox{}\linebreak[4]}

\newcommand{\N}{\mathbb N}
\newcommand{\Z}{\mathbb Z}
\newcommand{\Q}{\mathbb Q}
\newcommand{\R}{\mathbb R}
\newcommand{\C}{\mathbb C}

\newcommand{\MB}{{\mathcal B}}
\newcommand{\MC}{{\mathcal C}}
\newcommand{\MD}{{\mathcal D}}
\newcommand{\MF}{{\mathcal F}}

\newcommand{\MN}{{\mathcal N}}
\newcommand{\MNR}{{\mathcal R}}
\newcommand{\MS}{{\mathcal S}}
\newcommand{\MP}{{\mathcal P}}

\newcommand{\MO}{{\mathcal O}}
\newcommand{\SC}{\mathbf \Delta}
\newcommand{\GL}{\mathsf{GL}}
\newcommand{\gldz}{\GL_d(\Z)}   
  
\newcommand{\gldr}{\GL_d(\R)}   
\newcommand{\odr}{\mathsf{O}_d(\R)} 
\newcommand{\sd}{{\mathcal S}^d}
\newcommand{\sdo}{{\mathcal S}^d_{>0}}
\newcommand{\sdgeo}{{\mathcal S}^d_{\geq 0}}
\newcommand{\sdclosed}{\tilde{\mathcal S}^d_{\geq 0}}
\newcommand{\oK}{{\mathfrak o}_K}
\renewcommand{\vec}[1]{\boldsymbol{#1}}

\DeclareMathOperator{\trace}{trace}
\DeclareMathOperator{\sign}{sign}
\DeclareMathOperator{\vol}{vol}
\DeclareMathOperator{\conv}{conv}
\DeclareMathOperator{\cone}{cone}
\DeclareMathOperator{\aff}{aff}

\DeclareMathOperator{\BR}{BR}

\DeclareMathOperator{\Del}{Del}
\DeclareMathOperator{\vertex}{vert}
\DeclareMathOperator{\Min}{Min}

\DeclareMathOperator{\relint}{relint}
\DeclareMathOperator{\Aut}{Aut}
\DeclareMathOperator{\rank}{rank}

\theoremstyle{definition}
\newtheorem{definition}{Definition}[section]

\newtheorem{example}[definition]{Example}

\newtheorem{remark}[definition]{Remark}

\theoremstyle{plain}
\newtheorem{proposition}[definition]{Proposition}
\newtheorem{lemma}[definition]{Lemma}
\newtheorem{theorem}[definition]{Theorem}
\newtheorem{corollary}[definition]{Corollary}

\author{Mathieu Dutour Sikiri\'c}
\address{
Mathieu Dutour Sikiri\'c, Rudjer Boskovi\'c Institute, Bijenicka 54, 10000 Zagreb, Croatia}
\email{mdsikir@irb.hr}
\author{Achill Sch\"urmann}
\address{Achill Sch\"urmann, Mathematics Department, University of Magdeburg, 39106 Magdeburg, Germany}
\email{achill@math.uni-magdeburg.de}
\author{Frank Vallentin}
\address{Frank Vallentin, Centrum voor Wiskunde en Informatica (CWI), Kruislaan 413, 1098 SJ Amsterdam, The Netherlands} 
\email{f.vallentin@cwi.nl}
\thanks{The second and the third author were supported by the 
Deutsche Forschungsgemeinschaft (DFG) under grant SCHU 1503/4-1.
During the work on this paper the third author was partially supported
by the Edmund Landau Center for Research in Mathematical Analysis and
Related Areas, sponsored by the Minerva Foundation (Germany), and he was
partially supported by the Netherlands Organization for Scientific
Research under grant NWO 639.032.203.}

\subjclass{11H55, 52C17}

\title{A generalization of Voronoi's reduction theory and its application}

\begin{document}


\begin{abstract}
We consider Voronoi's reduction theory of positive definite quadratic
forms which is based on Delone subdivision.  We extend it to forms and
Delone subdivisions having a prescribed symmetry group. Even more
general, the theory is developed for forms which are restricted to a
linear subspace in the space of quadratic forms.  We apply the new
theory to complete the classification of totally real thin algebraic
number fields which was recently initiated by Bayer-Fluckiger and
Nebe.  Moreover, we apply it to construct new best known sphere
coverings in dimensions $9,\dots,15$.
\end{abstract}




\maketitle

\vspace*{-1.5cm}

\tableofcontents


\section{Introduction}

In this paper we generalize a classical reduction theory for positive
definite quadratic forms due to Voronoi \cite{voronoi-1908}.  His
theory gives in particular an algorithm to classify Delone subdivisions
of $\R^d$ with vertex-set $\Z^d$ up to the action of $\GL_d(\Z)$. For
precise definitions of used terms and a brief description of the
classical theory we refer to Section \ref{sec:background}.

We present our generalization in Section~\ref{sec:generalization}. We
extend the classical theory in two different directions. On the one
hand, we generalize the theory from vertex-set $\Z^d$ to general
periodic vertex-sets. On the other hand we give an equivariant theory
dealing with positive definite quadratic forms with a prescribed
automorphism group $G \leq \gldz$. This equivariant theory is an
analogue of the theory of $G$-perfect forms by A.M. Berg\'e,
J. Martinet and F. Sigrist, \cite{bms-1992}, which was motivated by
the search of good packing lattices with prescribed symmetries.  In
fact, as in their case, our theory can be developed in the more
general context of a linear subspace $T$ of quadratic forms.

As in the classical theory, where $T$ is the space of all quadratic
forms, we get a polyhedral subdivision of the space of positive
definite forms in $T$ into generic $T$-secondary cones (in the
classical case also called $L$-type domains), which contain those
forms giving the same Delone subdivision.  In contrast to the
classical theory these generic subdivisions are no longer
triangulations. In the equivariant theory, we have only finitely many
$T$-secondary cones up to the action of $\GL_d(\Z)$. Our proof of this
fact in Section~\ref{sec:equivariant} uses only the action of
$\GL_d(\Z)$ on a polyhedral subdivision of the space of positive
definite quadratic forms. So it applies to the theory of $G$-perfect
forms and it gives a unified view on both theories (see Remark
\ref{rem:unifiedview}).

We describe the theory in a way which allows us to work with it
computationally. In particular we made some effort to reduce
redundancies in the description of secondary cones (see Theorem
\ref{th:secondarycone}).  Also the transition from a Delone
subdivision of a $T$-secondary cone to the Delone subdivision of a
contiguous $T$-secondary cone, called a $T$-flip, is given explicitly
(see Theorem \ref{th:bistellar}).  For it we define repartitioning
polytopes, in which a polyhedral subdivision has to be replaced by
another. As a nice byproduct we obtain an explicit description of
flips which occur in the theory of equivariant secondary polytopes of
regular subdivisions (of a polytope), recently introduced by Reiner
\cite{reiner-2002} (see Remark \ref{rem:reiner}).

We are not the first who consider generalizations of Voronoi's
theory. Periodic tilings and reduction theory of positive definite
quadratic forms appear naturally in algebraic geometry in the study of
degenerations of abelian varieties and compactifications of Siegel
modular varieties. For these reasons, Delone subdivisions and
Voronoi's reduction theory were used and studied by Mumford and
Namikawa about 25 years ago, and by many algebraic geometers since
then. In Section \ref{sec:comparison} we review generalizations which
came up in this context and compare them to our work.

We apply our extension of Voronoi's theory to two different problems.
We use it to finish the classification of totally real thin number
fields, which was recently started by Nebe and Bayer--Fluckiger
\cite{bn-2004} (Section \ref{sec:thinfields}).  We use it to construct
new best known sphere coverings (Section \ref{sec:newcoverings}).
Both applications involve finding best lattice configurations with
respect to a given Delone subdivision. A brief description of the
problem and of the convex optimization tools we used is given in
Section~\ref{sec:coverings}. Since we are dealing with non-linear
optimization problems we can usually only approximate the lattices we
want to find. By using convex optimization duality and rational
approximations we can give mathematical rigorous error bounds for the
quality of approximated lattices. Algorithmic issues in the
classification of $T$-secondary cones are addressed in
Section~\ref{sec:algorithms}.

So far, using the new theory for the lattice case, we found new best
known lattice sphere coverings in dimensions $d=9,\dots,15$.  Using
the classical theory (cf. \cite{sv-2005}) and new methods to enumerate
all vertices of symmetric Voronoi cells efficiently
(cf. \cite{dsv-2006}), we found new best known coverings in dimensions
$d=6,7,8$ and $d=17,19,20,21$ as well.  A complete list of the best
known values is given in Table 2 in Section \ref{sec:newcoverings}.
With the exception of dimension $d=6$ and $d=7$, we do not think that
these lattice coverings are optimal.

Furthermore, we strongly believe that our extension of Voronoi's
classical theory to periodic sets is a first step towards the
construction of non-lattice coverings, which are less dense than any
lattice covering.


\section{Background: Lattices, PQFs and Delone subdivisions}

\label{sec:background}

We start with basic definitions and basic results and some background
on Voronoi's reduction theory.  In the first section we introduce
lattices and positive definite quadratic forms, PQFs from now on.  For
further reading we refer to \cite{cs-1988} and \cite{gl-1987}.  In the
second section we introduce Delone polyhedra and Delone
subdivisions. For further reading about concepts related to polyhedra
we refer to \cite{ziegler-1998}. In the third section we give a very
brief review on Voronoi's results in \cite{voronoi-1908}.

Let $\R^d$ be the $d$-dimensional Euclidean space with column vectors
$\vec{x} = (x_1, \ldots, x_d)^t$ and norm 
$\|\vec{x}\| = \sqrt{\vec{x}^t \vec{x}}$.

\subsection{Lattices and PQFs}

A $d$-dimensional \textit{lattice} $L$ in $\R^d$ is a discrete
subgroup $L = \Z\vec{v}_1 + \cdots + \Z\vec{v}_d$ with linearly
independent $\vec{v}_i\in\R^d$. The family
$(\vec{v}_1,\ldots,\vec{v}_d)$ is called a \textit{basis} of $L$. To
it we associate the positive definite symmetric Gram matrix $Q_B = B^t
B$, where $B\in\gldr$ is the invertible matrix whose $i$-th column is
$\vec{v}_i$, and $L=B\Z^d$.

Given on the other hand a positive definite symmetric matrix $Q$,
there exists a matrix $B\in\gldr$ with $Q=B^t B$.  The matrix $B$ is
uniquely determined up to orthogonal transformations.  Any other
$A\in\gldr$ with $A\Z^d=B\Z^d$ can be written as $A = B U$ with $U\in
\gldz$.  This relation yields $Q_{A} = A^t A = U^t Q_B U$.
We say $Q_{A}$ and $Q_B$ are \textit{arithmetical equivalent}
in this case. 

The space of real symmetric $d\times d$ matrices is denoted by $\sd$.
It is a $\binom{d+1}{2}$-dimensional Euclidean space with inner
product $\langle A, B \rangle = \trace(AB)$.  The subset $\sdo$ of
positive definite symmetric matrices is an open convex cone in $\sd$.
Abusing notation, we identify positive definite quadratic forms, PQFs
from now on, and positive definite symmetric matrices by
\[
Q[\vec{x}] = \vec{x}^t Q \vec{x} = \langle Q, \vec{x}\vec{x}^t\rangle
.
\]
The topological closure of $\sdo$ in $\sd$ is the set $\sdgeo$ of all
positive semidefinite matrices.  By the relations above, $\sdo$ can be
identified with $\odr \backslash \gldr$, where $\odr$ denotes the
subgroup of orthogonal $d\times d$ matrices in $\gldr$.  The group
$\gldz$ acts on $\sdo$ by $Q \mapsto U^t Q U$.  Thus the set of
\textit{isometry classes} of $d$-dimensional lattices (i.e.\ $d$-dimensional
lattices up to orthogonal transformations) can be identified with
$\sdo / \gldz$.

\subsection{Delone polyhedra and Delone subdivisions}

A \textit{polyhedron} is a set in $\R^d$ which can be represented as a
finite intersection of closed half spaces, e.g. given by a 
system of linear inequalities.
A \textit{polytope} is the convex hull of finitely many points
and by a theorem due to Minkowski and Weyl, polytopes are  
bounded polyhedra and vice versa. Given a discrete set $\Lambda\subset \R^d$, 
a polytope $P=\conv\{\vec{v}_1,\dots,\vec{v}_n\}$ with vertices 
$\vec{v}_i \in \Lambda$
is called a \textit{Delone polytope}, if there exists a center
$\vec{c}\in\R^d$ and a radius $r > 0$ such that the Euclidean distance
between $\vec{c}$ to all points $\vec{v}\in \Lambda$ satisfies
$\|\vec{c} - \vec{v}\|\geq r$, with equality only for the vertices of $P$.
The set of all Delone polytopes of a lattice forms a
\textit{polyhedral subdivision} of ${\R}^d$.
This is a family of polyhedra, called \textit{faces}, 
whose union is $\R^d$ and which is closed with 
respect to intersections. Each face of the subdivision sharing
relative interior points with another face of its dimension coincides
with this face. 
Note, by a theorem of Gruber and Ryshkov \cite{gr-1989},
the latter {\em face-to-face} property holds whenever
it holds for the facets (faces of co-dimension $1$).

Our main interest is in vertex-sets $\Lambda_L$ which are \textit{periodic},
that is, a finite union of lattice translates of a lattice $L$, e.g.,
$\Lambda_L=\bigcup_{i=1}^m \vec{t}_i + L$ with
$\vec{t}_i\in\R^d$ for $i=1,\dots,m$.
In many cases it is convenient to work with coordinates
with respect to a given basis of $L$, say given by $B\in\gldr$.
This means that we work with \textit{standard periodic sets} 
\begin{equation}
\label{eqn:standard_periodic_set} 
\Lambda=\bigcup_{i=1}^m \vec{t}'_i+\Z^d
\end{equation} 
and norm defined
by $Q_B$. Hence the norm of $\vec{x}\in\R^d$ is given by
$\sqrt{Q_B[\vec{x}]}$.
A polytope $P = \conv\{\vec{v}_1, \ldots , \vec{v}_n\}$ 
with $\vec{v}_i \in \Lambda$, is called a \textit{Delone polytope} of $Q_B$
and $\Lambda$, if there exists a $\vec{c}\in \R^d$ and a real number $r$ with 
$Q_B [\vec{c}-\vec{v}] \geq r^2$, where equality holds if and only if 
$\vec{v}$ is a vertex of $P$. 

It is a bit more general (and in some situations convenient)
to consider Delone subdivisions of positive semidefinite forms $Q$ as
proposed by Namikawa in \cite{namikawa-1976}. Positive semidefinite
forms define seminorms on $\R^d$ by $\sqrt{Q[\vec{x}]}$.  A (possibly
unbounded) polyhedron $P = \conv\{\vec{v}_1, \vec{v}_2, \ldots \}$
with $\vec{v}_i \in \Lambda$ is then called a \textit{Delone
polyhedron} of $Q$, if there exists a $\vec{c}\in \R^d$ and a real
number $r$ with $Q [\vec{c}-\vec{v}_i] = r^2$ for $i=1,2,\ldots$ and
$Q [\vec{c}-\vec{v}] > r^2$ for all $\vec{v}\in\Lambda\setminus
\{\vec{v}_1, \vec{v}_2, \ldots \}$.  Suppose for $Q\in\sdgeo$ there
exists a matrix $U\in\gldz$ and a $Q'\in{\MS}^{d'}_{>0}$ with
\begin{equation} 
\label{Q-reduce}
U^t Q U = \begin{pmatrix} 0 & 0 \\ 0 & Q'\end{pmatrix} = \bar{Q}
.\end{equation} Then $d'= \rank Q$ and if $P$ is a Delone polyhedron
of $Q$ then $U^{-1}P$ is a Delone polyhedron of $\bar{Q}$. Latter are
of the form $\R^{d-d'}\times P'$ with $P'\in\R^{d'}$ a Delone polytope
of $Q'$. The set of all forms being arithmetical equivalent to some
positive semidefinite form $\bar{Q}$ of the form \eqref{Q-reduce},
with $Q'$ positive definite, is called the \textit{rational closure}
of $\sdo$, denoted by $\sdclosed$. The importance of the rational
closure is due to the following proposition.

\begin{proposition}
\label{prop:sdclosed}
For $Q\in\sdgeo$ 
exists a $d$-dimensional Delone polyhedron
with respect to a standard periodic vertex-set
if and only if $Q\in\sdclosed$.
\end{proposition}

Although this proposition might be known we are not aware of a
reference. In \cite[\S 2.1]{namikawa-1976} Namikawa showed that every
$Q$ lying in the rational closure has a $d$-dimensional Delone
polyhedron. For completeness we repeat his argument.

\begin{proof}
If $Q\in\sdclosed$, there exist $d$-dimensional Delone polyhedra which
after a suitable transformation in $\GL_d(\Z)$ are those of the form
$\bar{Q}$ in \eqref{Q-reduce}.  If on the other hand
$Q\not\in\sdclosed$, then for all arithmetical equivalent forms
$\bar{Q}$ as in \eqref{Q-reduce}, the form $Q'\in \MS^{d'}_{\geq 0}$
is not positive definite, hence $d'>\rank(Q)$.  For an arithmetical
equivalent form $\bar{Q}$ with minimal $d'$ there exists no rational,
hence no integral vector in the kernel of $Q'$.  As a consequence, for
such a $Q'$ we find for all $\vec{c}\in\R^{d'}$ and all $r>0$ a
$\vec{v}\in\Z^{d'}$ with $Q'[\vec{v}-\vec{c}]<r$ (see
\cite[Lecture VI]{siegel-1989}).  Therefore there do not exist
$d'$-dimensional Delone polyhedra for $Q'$ and $\Z^{d'}$, and hence no
$d$-dimensional ones for $Q$ and $\Z^d$ as well.  The same is true for
$Q$ and a standard periodic set $\Lambda$.
\end{proof}

The set $\sdclosed$ can also be described as the set of all
non-negative combinations (the \textit{cone}) of rank-$1$ forms
$\vec{v}\vec{v}^t$ with $\vec{v}\in\Z^d$:

\begin{proposition}
\label{prop:rationalrays}
We have
\[
\sdclosed
= 
\cone \left\{ \vec{v}\vec{v}^t : \vec{v}\in\Z^d \right\}
.
\]
\end{proposition}

\begin{proof}
By the definition of $\sdclosed$, every $Q\in\sdclosed$ is
arithmetical equivalent to a form $\bar{Q}$ as in \eqref{Q-reduce}.
The PQF $Q'$ is a sum of rank-$1$ forms (cf. \cite[Section
24]{voronoi-1907}).  Therefore $\bar{Q}$ and $Q$ are of this form as
well.
 
Suppose on the other hand that $Q\in\sdgeo$ is the sum of rank
$1$-forms, e.g.  $Q=\sum_{i=1}^m \alpha_i \vec{v}_i \vec{v}_i^t$ with
$\vec{v}_i\in\Z^d$ and $\alpha_i>0$ for $i=1,\ldots,m$.  If $\rank Q =
d$ there is nothing to show.  If $\rank Q <d$, then there exist
linearly independent $\vec{p}_1,\ldots,\vec{p}_k\in\R^d$ with
$k=d-\rank Q$ such that $Q[\vec{p}_j]=\sum_{i=1}^m \alpha_i
(\vec{p}_j^t\vec{v}_i)^2 = 0$ for $j=1,\dots,k$. Thus the $\vec{p}_j$
are orthogonal to each of the $\vec{v}_i$ and therefore they span a
$k$-dimensional linear subspace which contains a $k$-dimensional
sublattice of $\Z^d$. We choose a basis $(\vec{u}_1,\ldots,\vec{u}_k)$
of this sublattice and extend it to a basis $U =
(\vec{u}_1,\ldots,\vec{u}_d)$ of $\Z^d$.  Then $(U^{-1})^t Q (U^{-1})$
is of the desired form \eqref{Q-reduce}.
\end{proof}

The set $\Del(Q)$ of all Delone polyhedra of a $Q\in\sdclosed$ is
called the \textit{Delone subdivision} of~$Q$.  If all elements in
$\Del(Q)$ are simplices, $\Del(Q)$ is called a
\textit{triangulation}. The subdivision $\Del(Q)$ is a polyhedral
subdivision of $\R^d$ which is invariant under translations of the
form $\vec{x}\mapsto\vec{x} + \vec{v}$, where $\vec{v}\in\Z^d$.
Therefore $\Del(Q)$ is completely determined by the \textit{stars} of
the translation vertices $\vec{t}'_i$ in
\eqref{eqn:standard_periodic_set} for $i=1,\dots,m$, where a star of a
single vertex is the set of all Delone polyhedra containing it.  We
call two Delone polyhedra $P$ and $P'$ \textit{equivalent} if there
exists a $\vec{v}\in\Z^d$ so that $P = \vec{v} + P'$. We say that
$\Del(Q)$ is a \textit{refinement} of $\Del(Q')$ (and $\Del(Q')$ is a
\textit{coarsening} of $\Del(Q)$), if every Delone polytope of $Q'$ is
contained in a Delone polytope of $Q$. In Section~\ref{sec:algorithms}
we need an algorithm which computes the Delone subdivision of a given
PQF. The interested reader can find a discussion of these
computational issues in our paper \cite{dsv-2006}.

\subsection{Voronoi's reduction theory}
\label{ssec:voros-theory}

Before we generalize Voronoi's reduction theory in the next section,
we briefly recall the original theory (see \cite{voronoi-1908}, \cite{delone-1937} and
\cite{sv-2005}).  Generally, the task of reduction is to find a
fundamental domain in $\sdo$ with respect to the action of $\gldz$.
Voronoi's reduction is based on \textit{secondary cones}, also called $L$-type
domains, of Delone triangulations with vertex-set $\Z^d$.  More generally, the
\textit{secondary cone} $\SC(\MD)$ of a Delone subdivision $\MD$ 
with a standard periodic vertex-set is defined by 
\[
\SC(\MD) = \{Q \in \sdclosed : \Del(Q) = \MD\}.
\]
We say that two secondary cones of Delone triangulations are
\textit{bistellar neighbors} if the Delone triangulations differ by a
\textit{bistellar flip}, which is a specific change of the triangulation
(see Section~\ref{ssec:bistellar} for a definition and generalization).
Voronoi also showed that the topological closures $\overline{\SC(\MD)}$ of
secondary cones of Delone triangulations form a polyhedral subdivision
of $\sdclosed$.

\begin{theorem}[\textit{Voronoi's Reduction Theory}] 
\label{th:mainvoronoi}
\finishline 
The secondary cone of a Delone triangulation with vertex-set $\Z^d$ 
is a full-dimensional, open polyhedral cone in $\sdo$.  The
topological closures $\overline{\SC(\MD)}$ give a polyhedral
subdivision of $\sdclosed$.  The closures of two secondary cones have
a common facet if and only if they are bistellar neighbors. The group
$\gldz$ acts on the tiling by $U \mapsto U^t \overline{\SC(\MD)}
U$. Under this group action there are only finitely many inequivalent
secondary cones.
\end{theorem}

Note that by Voronoi's theory we have a non-intersecting subdivision
of $\sdo$, as well as of $\sdclosed$, into secondary cones. In it,
every cone is an open polyhedral cone with respect to its affine 
hull. We refer to such a decomposition of $\sdo$ as an \textit{open
polyhedral subdivision}. Such subdivisions are of particular interest in
Section~\ref{sec:equivariant}, if they fall into only finitely many orbits under 
the action of $\gldz$, as in the case of secondary cones.


\section{Generalization of Voronoi's reduction theory}
\label{sec:generalization}

In this section we generalize Voronoi's reduction theory.
For our generalization we consider a linear subspace $T \subseteq
\MS^d$ and look at \textit{$T$-secondary cones} of Delone subdivisions
$\MD$ defined by
\[
\SC_T(\MD) = \SC(\MD) \cap T.
\]
We call a $T$-secondary cone $\SC_T(\MD)$ 
and the corresponding Delone subdivision $\MD$ 
\textit{$T$-generic} if $\dim \SC_T(\MD) = \dim T$.
By Voronoi's Theorem \ref{th:mainvoronoi}, the
topological closures of $T$-generic, $T$-secondary cones
give a polyhedral subdivision of $\sdclosed \cap T$. 
Two $T$-generic cones are called \textit{contiguous}
if their closures share a facet.
A difference with the classical theory is the existence 
of \textit{dead-ends}, which are facets only incident to
one $T$-generic cone. These necessarily contain only 
non-positive forms in $\sdclosed\setminus\sdo$.

The ultimate goal would be to state for
every subspace $T$ a theorem as Theorem \ref{th:mainvoronoi},
which deals with the case $T = \MS^d$. 
It turns out though that in general, this is not always possible.
If $\dim T = 1$, the intersection of $T$ with $\sdo$ contains a PQF $Q$ and all its 
multiples. In this case a generalized Theorem \ref{th:mainvoronoi}
is trivially true. Therefore, if not stated otherwise
we assume $\dim T\geq 2$ in what follows.

The ``road map'' for our generalization is the following: In
Section~\ref{ssec:secondary} and \ref{ssec:fundamental-proof} we
determine the secondary cone of an \textit{arbitrary} Delone
subdivision \textit{explicitly}, because in our more general setup we
have to deal with Delone subdivisions which are not Delone
triangulations.  Then in Section~\ref{ssec:bistellar} we generalize
the notion of bistellar neighbors to the new setting.  In the new
theory not every subspace $T$ gives a polyhedral subdivision of
$\sdclosed \cap T$ with only finitely many inequivalent $T$-secondary
cones (see Remark \ref{ex:infinite-many}).  For specific vertex-sets
and subspaces though, there exist only finitely many inequivalent
$T$-secondary cones.  Such a finiteness result is given in
Section~\ref{sec:equivariant} for $\Z^d$ and subspaces containing all
PQFs which are invariant under a given finite subgroup of $\gldz$.
Thus we obtain an equivariant version of Voronoi's reduction theory. 
Using modern terminology, our proofs not only generalize Voronoi's
theory, but also shorten his argumentation.

\subsection{Polyhedral description of secondary cones}
\label{ssec:secondary}

Let $\MD$ be a Delone subdivision with a standard periodic vertex-set $\Lambda$. 
In this section we want to
describe the secondary cone $\SC(\MD)$. It will turn out that
$\SC(\MD)$ forms a (relative) open polyhedral cone in $\sdclosed$. 
Theorem~\ref{th:secondarycone} gives the precise statement. This
result is an adaption of Voronoi's ``fundamental theorem'' 
(see \cite[\S 77]{voronoi-1908}),
which deals with the generic case of Delone triangulations 
with respect to the vertex-set $\Z^d$ (see also \cite[Section 5.1]{sv-2005}). Actually the first statement of the second part of Theorem~\ref{th:secondarycone} is not explicitly stated in Voronoi's work. It goes back to Nakamura \cite[Lemma 1.1]{nakamura-1975}.

We describe below the polyhedral cone explicitly 
by linear equalities and inequalities. The description is needed for our application
and therefore we put some effort into avoiding redundancies.
The linear equalities are coming from $d$-dimensional non-simplicial polyhedra in $\MD$. 
Hence, in the generic case of Delone triangulations there are no linear equalities. 
The linear inequalities are coming from $(d-1)$-dimensional polyhedra (facets) in $\MD$. 
For the formulation of these linear conditions we define for an affinely 
independent set $V\subseteq  \R^d$ of cardinality $d+1$ and a point
$\vec{w}\in\R^d$ the quadratic form
\begin{equation} \label{eqn:N-forms}
N_{V,\vec{w}} = \vec{w}\vec{w}^t - \sum_{\vec{v} \in V}
\alpha_{\vec{v}} \vec{v} \vec{v}^t
,
\end{equation}
where the coefficients $\alpha_{\vec{v}}$ are uniquely determined by
the affine dependency 
$\vec{w} = \sum_{\vec{v} \in V} \alpha_{\vec{v}} \vec{v}$ 
with $1 = \sum_{\vec{v} \in V} \alpha_{\vec{v}}$.

The following theorem generalizes 
Voronoi's ``fundamental theorem'' for Delone triangulations
to arbitrary polyhedral subdivisions. 
As vertex-sets $\Lambda\subset \R^d$ 
we allow standard periodic sets.
Moreover, we allow degenerate (unbounded) polyhedra. 
The \textit{vertex-set} of a polyhedron $P$, denoted by $\vertex P$,
is defined as the set $\Lambda\cap P$.

\begin{theorem}
\label{th:secondarycone}
\begin{enumerate}
\item[I.]
Let $\MD$ be a polyhedral subdivision of $\R^d$ with a
standard periodic vertex-set $\Lambda$.
Then the closure $\overline{\SC(\MD)}$ is a polyhedral cone in $\sdclosed$
and $\SC(\MD)$ is the set of all $Q\in\sdclosed$ satisfying
\begin{enumerate}
\item 
for every $d$-dimensional polyhedron $P\in\MD$ the equalities
\[
\langle N_{V,\vec{w}},Q \rangle = 0
,
\]
for one (which can be chosen arbitrarily) 
affinely independent set of $d+1$ vertices 
$V\subseteq \vertex P$ and all $\vec{w}\in \vertex P$;

\item
for every $(d-1)$-dimensional polyhedron $F\in \MD$ the inequality
\[
\langle N_{V\cup\{\vec{w}\},\vec{w}'},Q \rangle > 0
,
\]
for one (which can be chosen arbitrarily) 
affinely independent set of $d$ vertices $V\subseteq \vertex F$
and two vertices $\vec{w}\in \vertex P\setminus F$ and $\vec{w}'\in \vertex P'\setminus F$ 
of the two adjacent $d$-dimensional polyhedra $P,P'\in\MD$ with $F=P\cap P'$. 
\end{enumerate}
\item[II.]
The map $\MD \mapsto \overline{\SC(\MD)}$ gives an isomorphism between 
the poset of Delone subdivisions of $\Lambda$ ordered by coarsening
and the poset of closures of secondary cones ordered by inclusion.
The closures of all secondary cones of Delone subdivisions
form a polyhedral subdivision of $\sdclosed$.
\end{enumerate}
\end{theorem}

Note that different choices of $V$, $\vec{w}$ and $\vec{w}'$ 
for the inequalities could yield different conditions,
respectively different forms $N_{V\cup\{\vec{w}\},\vec{w}'}$. 
Nevertheless, these are the same 
on the linear subspace $U$ defined by the equalities. 
In other words, their orthogonal projections 
$\pi_U(N_{V\cup\{\vec{w}\},\vec{w}'})$ onto $U$ are all
positive multiples of a uniquely determined form 
$N_{\MD,F}\in U$ with $\langle N_{\MD,F},N_{\MD,F}\rangle = 1$.

Note also that for every $\vec{v} \in \R^d$ we have
\begin{equation} \label{eqn:translation-invariance}
N_{V+\vec{v},\vec{w}+\vec{v}} = N_{V,\vec{w}}
.
\end{equation}
Therefore, and because the vertex-set of the subdivision is assumed 
to be periodic,
the theorem gives only finitely many inequalities.
This shows that $\overline{\SC(\MD)}$ is a polyhedral cone.

Finally, let us remark that the theorem is valid
for arbitrary periodic sets, that is, finite unions of lattice
translates $\vec{t}_i+L$, if we replace $\sdclosed$ 
by $(A^{-1})^t\sdclosed (A^{-1})$ 
where $A\in\gldr$ defines the lattice $L=A\Z^d$.

\subsection{Proof of the fundamental theorem}
\label{ssec:fundamental-proof}

In this section we prove Theorem \ref{th:secondarycone}.
We first give two propositions which both deal with the redundancies
in the set of equations and inequalities, one would obtain
by considering all possible choices of $V$, $\vec{w}$ and $\vec{w}'$.
Proposition \ref{prop:delone-equalities} takes care of the equalities
and Proposition \ref{prop:reduce_description} of the inequalities.
The latter shows that the orthogonal projections 
of all the forms $N_{V\cup\{\vec{w}\},\vec{w}'}$ for a facet $F$,
onto the linear subspace defined by the equalities,
are unique up to positive multiples.

\begin{proposition} 
\label{prop:delone-equalities}
Let $Q \in \sd$ and $V\subset\R^d$ 
be an affinely independent set of cardinality $d+1$.  
Let $\vec{c} \in \R^d$ and $r >0$ be such that $Q[\vec{c} - \vec{v}] = r^2$ for all $\vec{v}\in V$.
Then
\[
Q[\vec{w} - \vec{c}] - r^2
= \langle Q, N_{V,\vec{w}} \rangle.
\]
\end{proposition}

\begin{proof}
The proof is straightforward. We have
\begin{eqnarray*}
Q[\vec{w} - \vec{c}] - r^2 & = & 
\langle Q, \vec{w}\vec{w}^t \rangle + \langle Q,
-2 \vec{w}\vec{c}^t + \vec{c}\vec{c}^t \rangle - r^2\\
& = &
\langle Q, \vec{w}\vec{w}^t \rangle + \sum_{\vec{v} \in V}
\alpha_{\vec{v}} \langle Q, -2\vec{v}\vec{c}^t + \vec{c}\vec{c}^t \rangle 
- r^2,
\end{eqnarray*}
with $\alpha_{\vec{v}}$ as in \eqref{eqn:N-forms}.
For each $\vec{v} \in V$ we use the equality $Q[\vec{v} -
\vec{c}] = r^2$ which is equivalent to $- \langle Q, \vec{v}\vec{v}^t
\rangle = \langle Q, -2 \vec{v}\vec{c}^t + \vec{c}\vec{c}^t \rangle -
r^2$. This yields the desired expression.
\end{proof}

\begin{remark}
\label{rem:chirotopes}
By Proposition \ref{prop:delone-equalities}
the sign of $\langle Q,  N_{V,\vec{w}}\rangle$
has the following interpretation: If it is positive,
then $\vec{w}$ lies outside the circumsphere of the points in $V$,
where the circumsphere is taken with respect to the norm 
induced by $Q$. If the sign is $0$, then $\vec{w}$ lies on the 
circumsphere, and if it is negative, then $\vec{w}$ lies inside
the circumsphere. In computational geometry this insphere/outsphere
test is conveniently formulated using oriented matroid 
terminology (cf. \cite[Chapter 1.8]{blswz-1999}):
Let $V=(\vec{v}_1,\dots,\vec{v}_{d+1})$ be affinely 
independent points in $\R^d$ with positive orientation
and let $\vec{w}\in\R^d$. Then the \textit{chirotope}
\[
\chi_{(\vec{v}_1,\ldots,\vec{v}_{d+1},\vec{w})}(Q)=
\sign 
\begin{vmatrix}
1 & \hdots & 1 & 1 \\ 
\vec{v}_1 & \hdots & \vec{v}_{d+1} & \vec{w} \\
Q[\vec{v}_1] & \hdots & Q[\vec{v}_{d+1}] & Q[\vec{w}]
\end{vmatrix}
\]
satisfies  
$\chi_{(\vec{v}_1,\ldots,\vec{v}_{d+1},\vec{w})}(Q)
=\sign\langle Q, N_{V,\vec{w}}\rangle.$
\end{remark}

\begin{remark}
\label{rem:regulators}
Voronoi's theory and his description of secondary cones of Delone triangulations
is based on linear forms $\varrho_{(L,L')}$ on $\sd$, 
called \textit{regulators}. 
Voronoi defines them for pairs of adjacent simplices $(L,L')$ sharing a facet 
in a Delone triangulation.
Let $\vec{w}'$ be the vertex of $L'$ which is not a vertex of $L$.
As in Proposition \ref{prop:delone-equalities}, let $V$ denote the vertex-set 
of $L$ and define $N_{V,\vec{w}'}$. Then Voronoi's regulator 
$\varrho_{(L,L')}$ is a positive multiple of 
$\langle N_{V,\vec{w}'}, \cdot \rangle$.
\end{remark}

\begin{proposition}
\label{prop:reduce_description}
Let $P$ be a $d$-dimensional polyhedron in $\R^d$. 
Let $U$ be the linear subspace of all $Q \in \sd$
satisfying $\langle N_{V,\vec{w}}, Q \rangle = 0$ for all affinely independent sets $V\subseteq \vertex P$ 
of cardinality $d+1$ and for all $\vec{w}\in \vertex P$. 
Let $F$ be a facet of $P$.
\begin{enumerate}
\item
Let $V$ and $V'$ be two sets of cardinality $d$ containing affinely independent vertices of $F$
and let $\vec{w}\in\vertex P \setminus F$ and $\vec{w}'\in\R^d$. Then
$$
\pi_U ( N_{V\cup\{\vec{w}\},\vec{w}'} )
=
\pi_U ( N_{V'\cup\{\vec{w}\},\vec{w}'} )
.
$$
\item
Let $V$ be a set of cardinality $d$ containing affinely independent vertices of $F$
and let $\vec{u},\vec{w}\in\vertex P \setminus F$ and $\vec{w}'\in\R^d$. Then 
$$
\pi_U ( N_{V\cup\{\vec{w}\},\vec{w}'} )
=
\pi_U ( N_{V\cup\{\vec{u}\},\vec{w}'} )
.
$$
\end{enumerate}
\end{proposition}

\begin{proof}
In both cases we will show that the difference of the two considered
forms lies in the orthogonal complement of $U$.
\begin{enumerate}
\item
Every pair of affinely independent sets $V, V' \subseteq \vertex F$ of
cardinality $d$ can be connected by a chain $V = V_1, \ldots, V_n =
V'$ of affinely independent sets $V_i \subseteq \vertex F$ of
cardinality $d$ such that $|V_i \cap V_{i+1}| = d-1$. So we can assume
$|V \cap V'| = d-1$. Setting $\alpha_{\vec{w}'}=1$ there exist unique
numbers $\alpha_{\vec{w}}$ and $\alpha_{\vec{v}}$ for $\vec{v} \in V$
defining a affine dependency between the points $\vec{w},
\vec{w}'$ and $\vec{v}$ in $V$. This defines the 
form $N_{V\cup\{\vec{w}\},\vec{w}'}$.

We define $\vec{v}_1$ by $\vec{v}_1 \in V \setminus V'$.
Since the affine hull of $V$ and of
$V'$ equals the affine hull of $F$ there exist numbers
$\beta_{\vec{v}}$ for $\vec{v} \in V'$ such that $\sum_{\vec{v} \in
V'} \beta_{\vec{v}} = 1$ and $\vec{v}_1 = \sum_{\vec{v} \in V'}
\beta_{\vec{v}} \vec{v}$. Thus we have an affine dependency
\[
\left\lbrace
\begin{array}{lcl}
0  & = & 
\displaystyle \alpha_{\vec{w}} + \alpha_{\vec{w}'} + \sum_{\vec{v} \in V \setminus
\{\vec{v}_1\}} \alpha_{\vec{v}} + \sum_{\vec{v} \in V'}
\alpha_{\vec{v}_1}\beta_{\vec{v}},\\
\vec{0} & = &  
\displaystyle \alpha_{\vec{w}} \vec{w} + \alpha_{\vec{w}'} \vec{w}' + \sum_{\vec{v}
\in V \setminus \{\vec{v}_1\}} \alpha_{\vec{v}} \vec{v} +
\sum_{\vec{v} \in V'} \alpha_{\vec{v}_1}\beta_{\vec{v}} \vec{v}
\end{array}
\right.
,
\]
which defines the form $N_{V'\cup\{\vec{w}\},\vec{w}'}$. This gives
\[
N_{V\cup\{\vec{w}\},\vec{w}'} - N_{V'\cup\{\vec{w}\},\vec{w}'} =
\alpha_{\vec{v}_1} \left(\vec{v}_1\vec{v}_1^t - \sum_{\vec{v} \in V'}
\beta_{\vec{v}} \vec{v}\vec{v}^t \right)
.
\]
Since $V\cup V'$ is a minimal affinely dependent set, we can choose an arbitrary 
vertex $\vec{w}$ of $P$ and find that the right hand side is a multiple of 
$N_{V'\cup\{\vec{w}\},\vec{v}_1}$.

\item 
We take a closer look at the difference
$N_{V\cup\{\vec{w}\},\vec{w}'} - N_{V\cup\{\vec{u}\},\vec{w}'}$
and show it is a multiple of $N_{V\cup\{\vec{w}\},\vec{u}}$.
Set $\alpha_{\vec{w}'}=1$ again and let $\alpha_{\vec{w}}$ and 
$\alpha_{\vec{v}}$ with $\vec{v} \in V$ be real numbers defining the
affine dependency between the points $\vec{w}, \vec{w}'$ and $\vec{v}$ in $V$.
This defines the form $N_{V\cup\{\vec{w}\},\vec{w}'}$. 
In the same way let $\alpha'_{\vec{w}'}=1$, 
$\alpha'_{\vec{u}}$ and $\alpha'_{\vec{v}}$ be real numbers defining 
$N_{V\cup\{\vec{u}\},\vec{w}'}$.

We set $\beta_{\vec{u}} = \alpha'_{\vec{u}}$, 
$\beta_{\vec{w}} = -\alpha_{\vec{w}}$  and $\beta_{\vec{v}} = \alpha'_{\vec{v}}-\alpha_{\vec{v}}$
for $\vec{v}\in V$. Then $\sum \beta_{\vec{v}} \vec{v} = \vec{0}$ and $\sum \beta_{\vec{v}}=0$
where the sums run through all $\vec{v}\in V\cup\{\vec{u},\vec{w}\}$.
Thus 
\begin{equation}
\label{eq:Nforms-difference}
N_{V\cup\{\vec{w}\},\vec{w}'} - N_{V\cup\{\vec{u}\},\vec{w}'} = \beta_{\vec{u}} N_{V\cup\{\vec{w}\},\vec{u}}.
\end{equation}
\end{enumerate}
\end{proof}

Note, if $\vec{u}$ and $\vec{w'}$ in the
last calculation lie in 
opposite halfspaces with respect to the affine plane through $V$, 
then $\beta_{\vec{u}}>0$ in \eqref{eq:Nforms-difference}.
Therefore, repeated application yields the following proposition,
which we use for the proof of Theorem \ref{th:secondarycone} and 
in Section \ref{ssec:bistellar} to prove Theorem \ref{th:bistellar}.

\begin{proposition}
\label{prop:regulator_main_formula}
Let $V_1,\ldots,V_m\subset\R^d$ be affinely independent sets 
of cardinality $d+1$ with $|V_i\cap V_{i+1}|=d$ for
$i=1,\ldots,m-1$. Let $\vec{w}\in\R^d$ and $V_i$ be on opposite sides
of $\aff(V_i\cap V_{i+1})$ for $i=1,\ldots,m-1$. Then
$$
N_{V_1,\vec{w}} = N_{V_m,\vec{w}} + \sum_{i=1}^{m-1} \alpha_i N_{V_i,\vec{v}_{i+1}}
$$
with $\vec{v}_{i+1}\in V_{i+1}\setminus V_i$ and positive constants $\alpha_i$, for $i=1,\ldots,m-1$.
\end{proposition}

With these propositions at hand, we can give a proof of the ``fundamental theorem''.

\begin{proof}[Proof of Theorem \ref{th:secondarycone}]
I.
We show that $\SC(\MD)$ is given by the set of 
listed linear equalities and inequalities.
By \eqref{eqn:translation-invariance} this implies that $\overline{\SC(\MD)}$ 
is a polyhedral cone, because the Delone subdivision 
induced by a $Q\in\sdclosed$ and a standard periodic vertex-set 
contains only finitely many Delone polyhedra up to $\Z^d$ invariant translations.

For $Q\in\SC(\MD)$ the linear equalities and inequalities are satisfied
by Proposition \ref{prop:delone-equalities}.

Conversely, let us assume $Q\in \sd$ satisfies the linear equalities
and inequalities for every polytope in $\MD$.  By Proposition
\ref{prop:delone-equalities} and Proposition
\ref{prop:reduce_description} we can assume that these are valid for
all possible choices of $V$, $\vec{v}$ and $\vec{w}$.  Note that for
the use of Proposition \ref{prop:reduce_description} it is crucial to
observe, that for two linear subspaces $U,U'$ of $\sd$ with $T=U\cap
U'$ we have $\pi_T = \pi_U\circ\pi_{U'} = \pi_{U'}\circ\pi_{U}$.

Let $P$ be a $d$-dimensional polyhedron in $\MD$. Let $V \subseteq
\vertex P$ be the vertex-set of a $d$-simplex.

We show that all possible inequalities $\langle N_{V,w}, Q \rangle
\geq 0$, where $\vec{w} \in \MD$, are implied by the inequalities
$\langle N_{V, \vec{v}}, Q \rangle \geq 0$, where $\vec{v} \in \MD$ is
either a vertex of $P$ or a vertex of an adjacent Delone polyhedron of
$P$. Assume $\vec{w}\in \MD\setminus P$. Then we choose a sequence of
adjacent $d$-simplices with vertex-sets $V_1,\dots,V_m$ in $\MD$
satisfying the requirements of Proposition
\ref{prop:regulator_main_formula}.  In addition we require that each
vertex-set is contained in the vertex-set of a fixed polyhedron of
$\MD$, in particular $V_1\subseteq \vertex P$ and $\vec{w}\in V_m$.
Note that this can be achieved by looking at a refining triangulation
of $\MD$, which we can choose arbitrarily.  By our assumption on $Q$,
we have $\langle N_{V_i,\vec{v}_{i+1}}, Q \rangle \geq 0$ with
equality if and only if $V_i,V_{i+1}$ are subsets of the same
polyhedron.  Thus by Proposition \ref{prop:regulator_main_formula} and
Proposition \ref{prop:delone-equalities} we see that the equation
$\langle N_{V_m,\vec{w}}, Q \rangle = 0$ implies $\langle
N_{V_1,\vec{w}}, Q \rangle > 0$.

Now we prove that a form $Q \in \sd$ satisfying all linear conditions
must be positive semidefinite. Let $\vec{w} \in \MD$ and $\lambda \in
\Z$ and write $\lambda \vec{w} = \sum_{\vec{v} \in V} \alpha_{\vec{v}}
\vec{v}$ with $1 = \sum_{\vec{v} \in V} \alpha_{\vec{v}}$. Then
\[
0 \leq \langle N_{V, \lambda \vec{w}}, Q \rangle = \lambda^2
Q[\vec{w}] - \sum_{\vec{v} \in V} \alpha_{\vec{v}} Q[\vec{v}],
\]
and hence $Q[\vec{w}] \geq 0$ since the coefficients
$\alpha_{\vec{v}}$ depend only affinely on $\lambda$.

Furthermore, $Q$ lies in the rational closure $\sdclosed$: The linear
conditions are rational and thus define a polyhedral cone in which
every extreme ray contains rational semidefinite quadratic forms. By
the argument in the proof of Proposition~\ref{prop:rationalrays} all
quadratic forms in such a polyhedral cone lie in $\sdclosed$.

We finally show that there exists a $\vec{c}\in\R^d$ and $r>0$ such
that $Q[\vec{w}-\vec{c}]\geq r^2$ with equality if and only if
$\vec{w}\in \vertex P$. A solution of the system of linear equations
$Q[\vec{v} - \vec{c}] = Q[\vec{v}' - \vec{c}]$, where $\vec{v},
\vec{v}' \in V$ with $V\subseteq \vertex P$, gives $\vec{c}$ and $r$. Applying
Proposition~\ref{prop:delone-equalities} shows that the inequalities 
$Q[\vec{w}-\vec{c}]\geq r^2$ are valid with equality if and only if
$\vec{w}\in\vertex P$.

II.  For the assertion on the isomorphism of posets given by the map
$\MD \mapsto \overline{\SC(\MD)}$, we need to verify that the Delone
subdivision $\MD$ is a true coarsening of $\MD'$ if and only if
$\overline{\SC(\MD)}$ is strictly contained in $\overline{\SC(\MD')}$.
This follows from I. because if $\MD$ is a true coarsening of $\MD'$
we have equalities in the description of ${\SC(\MD)}$ which are
inequalities in the description of ${\SC(\MD')}$.  If on the other
hand, $\overline{\SC(\MD)}$ is strictly contained in
$\overline{\SC(\MD')}$, we know that all equalities and inequalities
in the description of $\SC(\MD')$ are also satisfied by elements of
$\SC(\MD)$ implying that $\MD$ is a coarsening of $\MD'$.  Moreover,
$\overline{\SC(\MD)}$ has to be contained in the boundary of
$\overline{\SC(\MD')}$, because otherwise there would exist a $Q\in
\SC(\MD')\cap \SC(\MD)$ implying $\MD=\MD'$.  Thus at least one of the
inequalities in the description of $\overline{\SC(\MD')}$ is fulfilled
with equality for the elements of $\overline{\SC(\MD)}$.  Hence $\MD$
is a \textit{true} coarsening of $\MD'$.

Since every form in $\sdclosed$ defines a unique Delone subdivision of
$\Lambda$, $\sdclosed$ is subdivided into secondary cones.  This
subdivision has the face-to-face property, hence is a polyhedral
subdivision, due to the poset isomorphism: Otherwise we would find
Delone subdivisions $\MD$ and $\MD'$, such that $\overline{\SC(\MD)}$
and $\overline{\SC(\MD')}$ have the same dimension and such that there
is a common relative interior point, but with
$\overline{\SC(\MD)}\not=\overline{\SC(\MD')}$.  Thus $Q$ would yield
two different Delone subdivisions, which is not possible.
\end{proof}

\subsection{Flips and bistellar neighbors}
\label{ssec:bistellar}

Given a linear subspace $T$ of $\sd$, we know by 
Theorem~\ref{th:secondarycone} that $\sdclosed\cap T$
is covered by the topological closures of 
$T$-generic $T$-secondary cones. 
As before 
we assume that the vertex-set of the considered Delone subdivisions
is a fixed standard periodic set $\Lambda$.
Given a $T$-generic Delone subdivision $\MD$ and its
$T$-secondary cone $\SC_T(\MD)$ we want to determine its
contiguous $T$-secondary cones. 
For this let $\MF$ be a facet of $\overline{\SC_T(\MD)}$ which is not a dead-end.
Let $\MD'$ be the $T$-generic Delone subdivision
with $\SC(\MD')$ being contiguous to $\SC(\MD)$ at $\MF$, hence with 
$\MF= \overline{\SC_T(\MD)}\cap \overline{\SC_T(\MD')}$. 
The transition from $\MD$ to $\MD'$ is called a 
\textit{$T$-flip}; the subdivisions $\MD$ and $\MD'$ are referred to as  
\textit{bistellar neighbors}.

In order to describe the $T$-flip, let $N_{\MF}\in\sd$ denote 
the form which is uniquely determined by the following conditions:  
\begin{enumerate}
\item[(i)]
$N_{\MF} \in T$ with $\langle N_{\MF},N_{\MF}\rangle=1$, 
\item[(ii)] 
$\langle N_{\MF}, Q \rangle = 0$ for all $Q \in \MF$, 
\item[(iii)] 
$\langle N_{\MF}, Q \rangle > 0$ for all $Q \in \SC_T(\MD)$, 
\item[(iv)] 
$\langle N_{\MF}, Q \rangle < 0$ for all $Q \in \SC_T(\MD')$. 
\end{enumerate}
We collect facets $F$ of $\MD$ whose forms $N_{\MD, F}$ 
(as defined after Theorem \ref{th:secondarycone}) are
positive multiples of $N_{\MF}$ when projected onto $T$:
\[
\MNR_{\MF} = \{F \in \MD : \mbox{$\dim F = d - 1$ and $N_{\MF} = 
\alpha \cdot \pi_T(N_{\MD, F})$ for a $\alpha > 0$}\}.
\]
On $\MNR_{\MF}$ we define an equivalence relation by
\begin{equation} 
\label{eqn:equivalence-rel}
F \sim F' \Longleftrightarrow
\left\{
\begin{array}{l}
\mbox{there exist $F_1, \ldots,  F_n \in \MNR_{\MF}$ with $F = F_1$, $F' = F_n$,}\\
\mbox{and $d$-dimensional polyhedra $P_0, \ldots, P_n \in \MD$}\\
\mbox{with $F_i = P_{i-1} \cap P_i$ for
$i = 1, \ldots, n$.}\\
\end{array}
\right.
\end{equation}

Thus by definition, the facets in each equivalence class are connected by a chain
of adjacent $d$-dimensional polyhedra in $\MD$. The union of these polyhedra
is a polyhedron again:

\begin{proposition}
\label{prop:repartitioning_polytopes}
Let $\MD$ be a Delone subdivision
and $\MC$ an equivalence class of \eqref{eqn:equivalence-rel}. 
Let $V_{\MC}$ be the set of all vertices of $d$-dimensional 
polytopes in $\MD$ with a facet in $\MC$. 
Then $(\conv V_{\MC})\cap \Lambda= V_\MC$.
\end{proposition}

A proof of the proposition is given at the end of this section.
It, as well as the description of $T$-flips
can conveniently be given using the \textit{lifting map}
$w_Q : \R^d \to \R^{d + 1}$ defined by $w_Q(\vec{x}) = (\vec{x}, Q[\vec{x}])$
for a $Q\in\sdclosed$.
The polyhedra $\conv V_{\MC}$ in Proposition \ref{prop:repartitioning_polytopes} 
are called \textit{repartitioning polyhedra} of $\MF$. 
For a $Q \in \SC_T(\MD)$ and a repartitioning polyhedron $P$
we define the $(d+1)$-dimensional polyhedron 
$$
\bar{P}(Q) = 
\conv\{w_Q(\vec{v}) : \vec{v} \in \vertex P\}.
$$
For each of its facets $\bar{F}$
we have an outer normal vector $\vec{n}(\bar{F}) \in \R^{d+1}$.
The facets are divided into three groups. We speak of a
\textit{lower facet}, if the last coordinate of the normal vector satisfies 
$(\vec{n}(\bar{F}))_{d+1}<0$; we speak of 
an \textit{upper facet}, if $(\vec{n}(\bar{F}))_{d+1}>0$
and of a \textit{lateral facet}, if $(\vec{n}(\bar{F}))_{d+1}=0$.
We show below that these notions are independent of the particular choice of
$Q\in \SC_T(\MD)$.

For the description of the $T$-flip in the following theorem,
let $\pi : \R^{d+1}\mapsto \R^d$ denote the projection
onto the first $d$ coordinates. Figure 1. gives an example 
for the change from ``lower to upper hull'' in the case of a planar repartitioning polytope. 
Note that the situation in higher dimensions can be much more complicated than
the picture might suggest.

\bigskip
\begin{center}
\includegraphics[width=3.5cm]{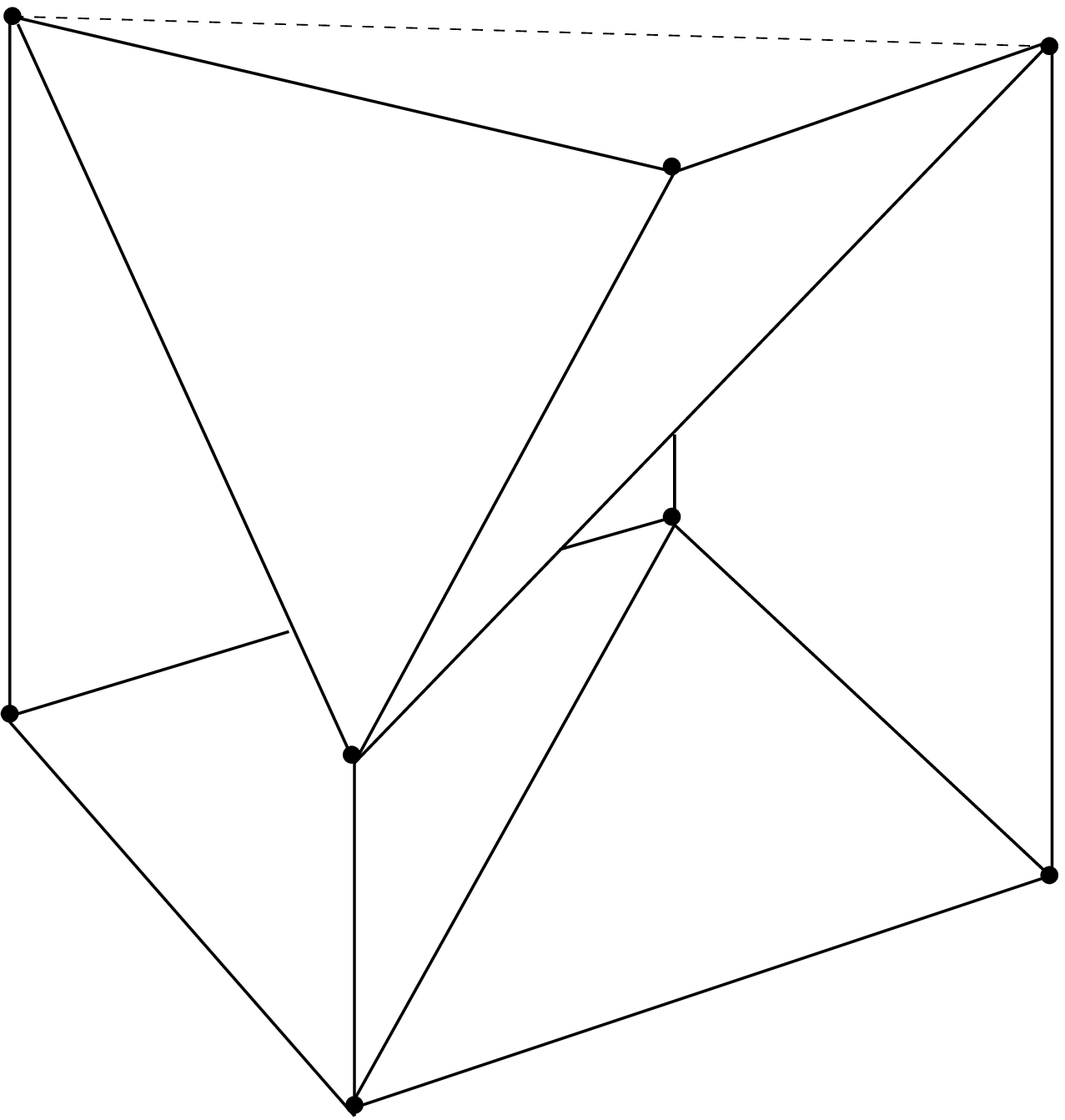}   
\hspace{2.0cm}
\includegraphics[width=3.5cm]{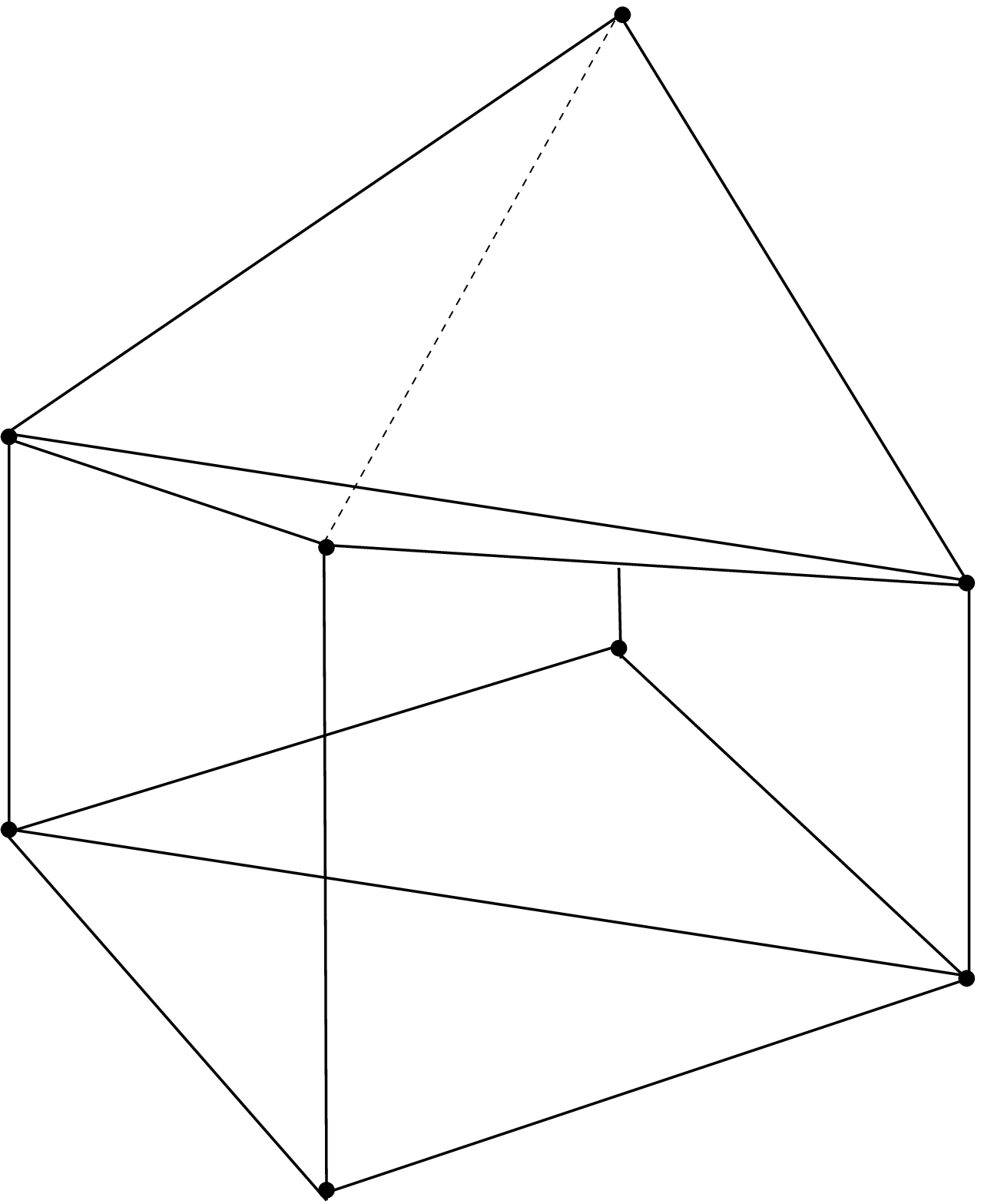}
\par\smallskip
\textbf{\textsf{Figure 1.}} Lifting interpretation of ``upper to lower hull change''.
\end{center}

\begin{theorem}
\label{th:bistellar}
Let $T$ be a linear subspace of $\sd$ and 
$\MF$ be a common facet of two contiguous $T$-secondary cones
of $T$-generic Delone subdivisions $\MD$ and $\MD'$. Then
we obtain the $d$-dimensional polytopes of $\MD'$ from those of $\MD$ by
choosing an arbitrary $Q\in\SC(\MD)$ and 
by doing the following for each repartitioning polytope $P$ of $\MF$:
\begin{enumerate}
\item[(i)]
We remove all polytopes $\pi(\bar{F})$ 
of lower facets $\bar{F}$ of $\bar{P}(Q)$.

\item[(ii)]
We add all polytopes $\pi(\bar{F})$ 
of upper facets $\bar{F}$ of $\bar{P}(Q)$.
\end{enumerate}
\end{theorem}

\begin{remark}
\label{rem:reiner}
Our theory is similar to the 
theory of equivariant secondary polytopes
of regular polytopal subdivisions, 
as recently described by Reiner in \cite{reiner-2002}. 
With slight modifications our explicit description 
of the flipping procedure works in his setting as well.
The linear subspace $T$ is in this case 
containing all PQFs which are invariant under a given 
finite subgroup of $\gldz$
(see Section \ref{sec:equivariant}).
\end{remark}

As a preparation for the proofs of Theorem \ref{th:bistellar}
and Proposition \ref{prop:repartitioning_polytopes}, 
consider the set
\begin{equation} \label{eqn:lifting_polyhedron}
\conv\{w_Q(\vec{x}) : \vec{x}\in \Lambda \}
.
\end{equation}
It is a \textit{locally finite polyhedron}, meaning that the intersection with
a polytope (bounded polyhedron) is a polytope again.
A first and important observation is that the facets of the
polyhedron \eqref{eqn:lifting_polyhedron}
yield the $d$-polytopes of the Delone subdivision when
projected by $\pi$ onto the first $d$ coordinates (see Figure 2):

\begin{proposition}
\label{prop:lifting_interpret}
Let $\Lambda$ be a standard periodic set and $Q\in\sdclosed$. Then
\[
\Del(Q)
=
\{\pi(\bar{F}) : \bar{F} \mbox{ face of \eqref{eqn:lifting_polyhedron} of dimension less than or equal $d$}\}
.
\]
\end{proposition}

\bigskip
\begin{center}
\includegraphics[width=8cm]{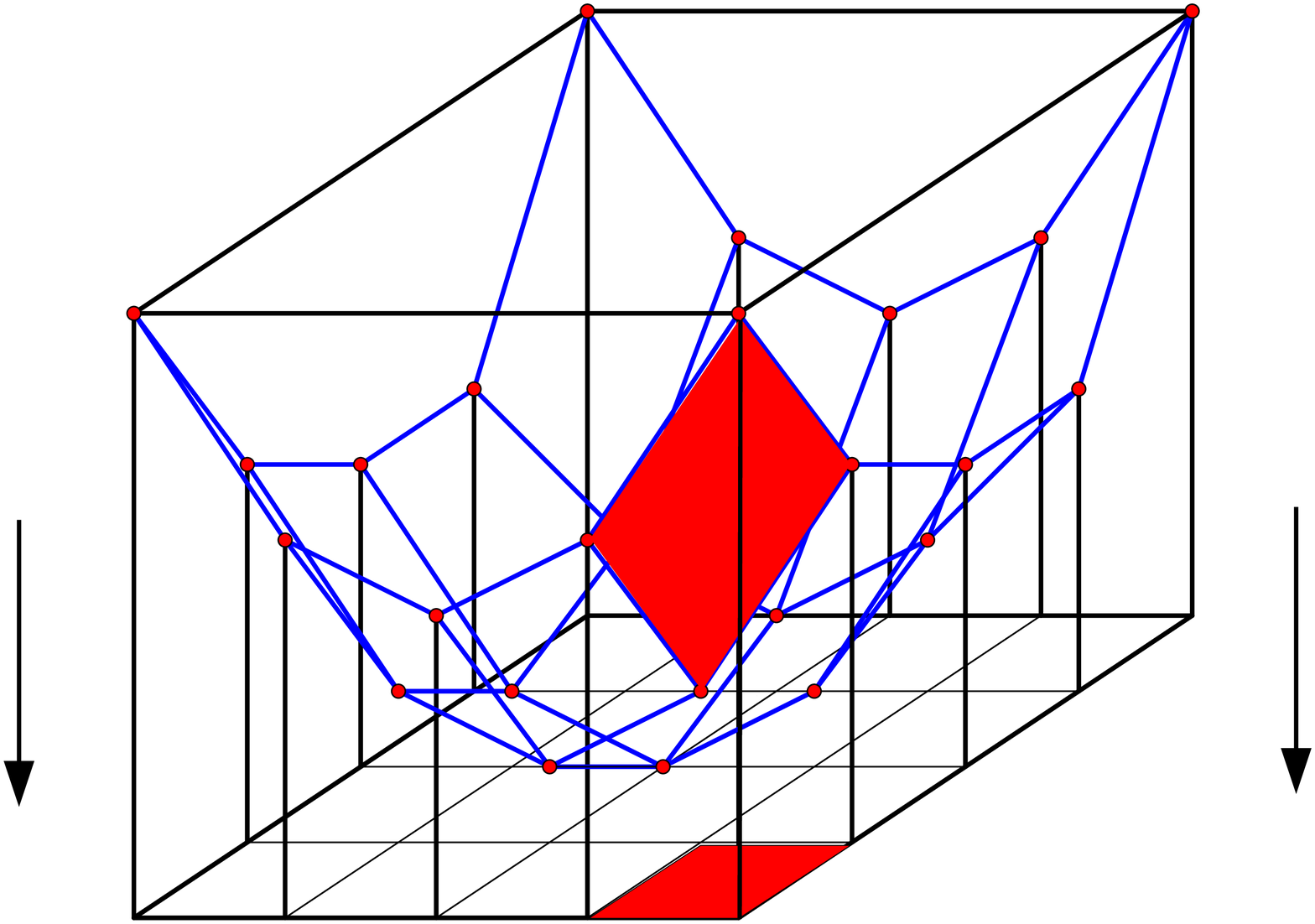}
\par\smallskip
\textbf{\textsf{Figure 2.}} Delone subdivision obtained from lifting.
\end{center}

\begin{proof}
For a polytope $P\in\MD$ and a set 
$V\subseteq \vertex P$ of cardinality $d+1$ and with affinely independent vertices,
the set of all $(\vec{x},y)\in\R^{d+1}$ with 
\begin{equation} \label{eqn:support_plane}
y=Q[\vec{x}]- \langle N_{V,\vec{x}},Q \rangle
 =\sum_{\vec{v}\in V} \alpha_{\vec{v}} Q[\vec{v}]
\end{equation}
is an affine hyperplane in $\R^{d+1}$.
Recall from the definition of $N_{V,\vec{x}}$ in  
\eqref{eqn:N-forms} that the $\alpha_{\vec{v}}$ 
are uniquely defined by the affine dependency 
$\vec{x}=\sum_{\vec{v}\in V} \alpha_{\vec{v}} \vec{v}$ with $\sum_{\vec{v}\in V} \alpha_{\vec{v}}=1$.
Thus we see that the hyperplane given by \eqref{eqn:support_plane}
contains the $d+1$ lifted points $w_Q(\vec{v})$.
Moreover, it is a supporting hyperplane of the polyhedron 
\eqref{eqn:lifting_polyhedron} if and only if 
$\langle N_{V,\vec{x}},Q \rangle\geq 0$ for all $\vec{x}\in\Lambda$.
By Proposition \ref{prop:delone-equalities} this is the case
if and only if $Q[\vec{x}-\vec{c}]\geq r^2$ for a 
suitable $\vec{c}\in\R^d$ and $r>0$. Here, 
equality holds if and only if $\vec{x}\in\vertex P$.
\end{proof}

The hyperplane given by equation \eqref{eqn:support_plane} in the proof
is a different one for subsets $V$ and $V'$ of different Delone polytopes.
The forms $N_{\MD,F}$ relate the change 
for adjacent $d$-polytopes $P$ and $P'$ in $\MD$ with common facet $F=P\cap P'$. 
If the segment connecting $\vec{x}$ with a vertex $\vec{v}$ of $P$
intersects facets $F_1,\dots,F_n$ of $\MD$, then
by Proposition \ref{prop:regulator_main_formula} 
and by the definition of $N_{\MD,F}$
(after and because of Theorem \ref{th:secondarycone}) 
we have
\begin{equation} \label{eqn:regulator_relation}
y=
Q[\vec{x}] - \langle N_{V,\vec{x}},Q \rangle
=
Q[\vec{x}] - \sum_{i=1}^n \beta_i \langle N_{\MD,F_i},Q \rangle
\end{equation}
with suitable positive constants $\beta_i$,
depending on $\vec{x}$.

Equation \eqref{eqn:regulator_relation} yields simple proofs for 
Proposition \ref{prop:repartitioning_polytopes} and Theorem \ref{th:bistellar}:

\begin{proof}[Proof of Proposition \ref{prop:repartitioning_polytopes}]
Due to \eqref{eqn:regulator_relation} and the definition
of $\MC$, all points $w_Q(\vec{v})$ with $\vec{v}\in V_\MC$
are coplanar in $\R^{d+1}$. Suppose the union of the $d$-dimensional 
polyhedra of $\MD$ with a facet in $\MC$ is strictly contained in 
$\conv V_\MC$. Then there exists a facet $F$ of one of the polyhedra
such that $F\not\in\MC$ and $F$ does not belong to the boundary
of $\conv V_\MC$. Hence, there exist two vertices $\vec{w},\vec{w}'\in V_\MC$
in opposite halfspaces with respect to $\aff F$.
Then by \eqref{eqn:regulator_relation} and since 
$\langle N_{\MD,F}, Q \rangle > 0$ we see that 
$w_Q(\vec{w})$, $w_Q(\vec{w}')$ and $w_Q(\vec{v})$ with $\vec{v}\in F$
can not be coplanar, which is a contradiction.
\end{proof}

\begin{proof}[Proof of Theorem \ref{th:bistellar}]
Choose $Q\in \SC(\MD)$ and $Q'\in\SC(\MD')$ and let
$P$ be a repartitioning polytope of $\MF=\overline{\SC(\MD)}\cap\overline{\SC(\MD')}$.
Then by formula \eqref{eqn:regulator_relation}
and since $N_{\MD,F} = - N_{\MD',F}$ for all facets $F$ contained in $P$,
the upper facets of $\bar{P}(Q)$ project to the same polyhedral 
subdivision of $P$ as the lower facets of $\bar{P}(Q')$
and vice versa. The lower facets give the corresponding Delone
subdivisions by Proposition \ref{prop:lifting_interpret}, which proves the assertion.
\end{proof}


\section{Finiteness and equivariance}
\label{sec:equivariant}

For the sake of simplicity and because it is sufficient for our applications
in Section \ref{sec:thinfields} and Section \ref{sec:newcoverings}
we restrict ourselves from now on to the case of Delone subdivisions
with vertex-set $\Z^d$.
A more general discussion can be found in \cite{schuermann-2007}.
In particular it is shown that the results of this Section extend
to the case of \textit{rational standard periodic} vertex-sets
$\Lambda=\bigcup_{i=1}^m \vec{t}_i + \Z^d$ with $\vec{t}_i\in\Q^d$
for $i=1,\dots,m$.

Let $T$ be a linear subspace of $\MS^d$. We say that two $T$-secondary
cones $\Delta_T$ and $\Delta'_T$ are \textit{T-equivalent} if there is
a $g \in \GL_d(\Z)$ so that $g^t \Delta_T g = \Delta'_T$ and $g^tTg=T$. 
Otherwise we say that they are \textit{T-inequivalent}. 
In this section we discuss
assumptions on~$T$ which ensure that there exist only finitely many,
$T$-inequivalent $T$-secondary cones. Note that Voronoi's classical
theory (cf. Theorem~\ref{th:mainvoronoi}) deals with the case $T = \sd$.

The main difference with the classical theory is that for general $T$
there may be infinitely many, $T$-inequivalent $T$-secondary cones:

\begin{example}
\label{ex:infinite-many}
Consider a rational PQF $Q$ with trivial automorphism group. The orbit
of $Q$ under the action of $\gldz$ contains only rational PQFs. We
choose a subspace $T$ through $Q$ in which the only rational forms are
multiples of $Q$. So $T$ does not contain another element from the orbit of $Q$. 
Thus the setwise stabilizer of $T$ is trivial. 
 In this situation there are infinitely many
$T$-generic secondary cones whenever the set $T \cap \sdclosed$ is not
closed in $\MS^d$. Then the intersection of $T$ with
the boundary of $\sdclosed$ is not covered by finitely many $T$-dead-ends.  
The maximum dimension of the
intersection of $\sdgeo$ with a supporting hyperplane of $\sdgeo$ is
$\binom{d+1}{2}-d$ (see \cite[\S 10]{br-1979}). Thus the intersection
of $T$ with the boundary of $\sdclosed$ cannot be covered by $T$-dead-ends 
if $\dim(T \cap \sdgeo) > \binom{d+1}{2}-d+1$. 
In particular for all $d \geq 3$ there exist
examples of subspaces $T$ with infinitely many, inequivalent
$T$-generic cones.
\end{example}


In the following we show that there are only finitely many,
$T$-inequivalent $T$-secondary cones if $T$ is the linear subspace which is
stabilized pointwise by a finite subgroup $G$ of $\GL_d(\Z)$. We
develop this equivariant theory in analogy to the theory of
$G$-perfect forms of Berg\'e, Martinet and Sigrist
\cite{bms-1992}. In Remark~\ref{rem:unifiedview} we give a unifying
view on both theories.

First let us recall some definitions and basic results (see
\cite[Section 2]{bnz-1973}). Let $G$ be a finite subgroup of
$\GL_d(\Z)$. The linear subspace
\[
\MF(G) = \{Q \in \MS^d : \mbox{$g^t Q g = Q$ for all $g \in G$}\}
\]
is called the \textit{space of invariant forms} of $G$. The pointwise
stabilizer of $\MF(G)$, which we denote by
\[
\MB(G) = \{g \in \GL_d(\Z) : \mbox{$g^t Q g = Q$ for all $Q \in \MF(G)$}\},
\]
is called the \textit{Bravais group} of $G$.  Note that $\MB(G)$ can
be strictly larger than $G$. The \textit{normalizer} of a subgroup
$G$ of $\GL_d(\Z)$ is defined by $N(G) = \{n \in \GL_d(\Z) : n^{-1} G
n = G\}$. One important property of $N(\MB(G))$ is that it is the
setwise stabilizer of the subspace of invariant forms
\[
N(\MB(G)) = \{g \in \GL_d(\Z) : g^t \MF(G) g = \MF(G)\}.
\]
A proof can be found for example in \cite[Lemme 3.2]{jaquet-1995}.

Now we are ready to state the main result of this section.  For this
recall that an open polyhedral subdivision is a non-intersecting 
decomposition into polyhedral cones which are open
with respect to their affine hull.

\begin{theorem}
\label{th:equivariance}
Let $\MP$ be an open polyhedral subdivision of $\sdo$ on which 
$\GL_d(\Z)$ acts by $(g, \Delta) \mapsto g^t \Delta g$. Suppose 
that this action gives only finitely many orbits. Define a polyhedral
subdivision of $\MF(G) \cap \sdo$ by
\[
\MP_{\MF(G)} = \{\Delta \cap \MF(G) : \Delta \in \MP\}.
\]
Then the normalizer $N(\MB(G))$ acts on $\MP_{\MF(G)}$ and there exist only 
finitely many orbits with respect to this action.
\end{theorem}

This together with Theorem~\ref{th:mainvoronoi} gives the following
corollary, which completes our equivariant version of Voronoi's theory.

\begin{corollary}
\label{cor:equivariance}
Let $G$ be a finite subgroup of $\GL_d(\Z)$ and let $T = \MF(G)$ be
the space of invariant forms. Then there exist only finitely many
$T$-inequivalent $T$-secondary cones of Delone subdivisions of $\Z^d$.
\end{corollary}

\begin{proof}[Proof of Theorem~\ref{th:equivariance}]
We start with a definition. For $\Delta \in \MP$ we define the
automorphism group $\Aut(\Delta) = \{g \in \gldz: g^t \Delta g =
\Delta\}$.

We consider the set
\[
\MO_{\Delta, G} = \{g^t \Delta g \cap \MF(G) : \mbox{$g \in \GL_d(\Z)$
and $g^t \Delta g \cap \MF(G) \neq \emptyset$}\}.
\]
The normalizer $N(\MB(G))$ stabilizes $\MF(G)$ setwise, and hence it
acts on $\MO_{\Delta, G}$. We show that $\MO_{\Delta, G}$ is a finite
union of $N(\MB(G))$-orbits, that is, there are $g_1, \ldots, g_m \in
\GL_d(\Z)$ with
\[
\MO_{\Delta, G} = \bigcup_{i = 1}^m
\{ n^t (g_i^t \Delta g_i \cap \MF(G)) n: n \in N(\MB(G)) \}.
\]
Then the statement of the theorem follows because the set
$\{\MO_{\Delta,G} : \Delta \in \MP\}$ is finite by the assumption made
on $\MP$.

Let $g_1, g_2$ be in $\GL_d(\Z)$ so that $g_i^t \Delta g_i \cap \MF(G)
\in \MO_{\Delta, G}$ for $i = 1,2$. The group $\MB(G)$ is a subgroup
of $\Aut(g_1^t \Delta g_1)$ and $\Aut(g_2^t \Delta g_2)$ which can be
seen as follows. By definition the group $\MB(G)$ stabilizes the set
$g_i^t \Delta g_i \cap \MF(G)$ pointwise. Since $\GL_d(\Z)$ operates
on $\MP$ and $g_i^t \Delta g_i \cap \MF(G) \neq \emptyset$ every
element of $\MB(G)$ has to stabilize $g_i^t \Delta g_i$ setwise.

Hence, $g_1 \MB(G) g_1^{-1}$ and $g_2 \MB(G) g_2^{-1}$ are subgroups of
$\Aut(\Delta)$. Assume that $g_1 \MB(G) g_1^{-1} = g_2 \MB(G) g_2^{-1}$. Then $n
= g_2^{-1}g_1$ is an element of the normalizer
$N(\MB(G))$. Furthermore, we have
\[
n^t (g_2^t \Delta g_2 \cap \MF(G)) n = n^t g_2^t \Delta g_2 n \cap n^t
\MF(G) n = g_1^t \Delta g_1 \cap \MF(G)
\]
because $N(\MB(G))$ stabilizes $\MF(G)$ setwise. Thus different
$N(\MB(G))$-orbits in $\MO_{\Delta, G}$ induce different subgroups of
$\Aut(\Delta)$. 

For every $\Delta \in \MP$ the group $\Aut(\Delta)$ is finite
(This fact was already proved by Nakamura in \cite[Lemma
1.2]{nakamura-1975}. For completeness we give an argument, which
also is a bit more elementary than Nakamura's.): Let $\MNR_1,\ldots, \MNR_k \in
\sdclosed$ be the extreme rays spanning the closed polyhedral cone
$\overline{\Delta}$. Choose forms $Q_i \in \MNR_i$ in the following
way: If for a pair $\MNR_i$, $\MNR_j$ there is a $U \in
\Aut(\Delta)$ with $U^t\MNR_iU = \MNR_j$ then $U^t Q_i U = Q_j$. We
have $\overline{\Delta}= \cone\{Q_1,\ldots, Q_k\}$ and $Q=\sum_{i=1}^k
Q_i \in \Delta$ and in particular $Q \in \sdo$.  So, $\Aut(Q) = \{g
\in \GL_d(\Z) : g^t Q g = Q\}$ is finite. By construction $\Aut(Q)$
contains $\Aut(\Delta)$. So, the group $\Aut(\Delta)$ has only a
finite number of different subgroups. So, $\MO_{\Delta, G}$ is a
finite union of $N(\MB(G))$-orbits.
\end{proof}

We close this section with two remarks.

\begin{remark}
\label{rem:unifiedview}
The proof of Theorem~\ref{th:equivariance} is essentially an
adaptation of proofs of Jaquet-Chiffelle
(\cite[Th\'eor\`eme~5.2]{jaquet-1995}) and Opgenorth (\cite[Theorem
4.2.18]{opgenorth-1995}).  They show that there exist only finitely
many inequivalent \textit{$G$-perfect forms} (definition is given
below).  Proofs for their theorems can be derived from
Theorem~\ref{th:equivariance} in the following way:

Let $m$ be a positive number and let $\MP_m$ be the set
\[
\MP_m = \{Q \in \sdo : \mbox{$Q[\vec{v}] \geq m$ for all $\vec{v} \in
\Z^d \setminus \{\vec{0}\}$}\}.
\]
This is a convex, locally finite polyhedral cone. For each face $F$ of
$\MP_m$ we define the (relatively) open polyhedral cone $\Delta_F =
\{\lambda Q : \mbox{$\lambda > 0$ and $Q \in \relint F$}\}$. Then $\MP
= \{\Delta_F : \mbox{$F$ face of $\MP_m$}\}$ gives an open polyhedral
subdivision of $\sdo$ as required in Theorem~\ref{th:equivariance}.

If $G$ is a finite subgroup of $\GL_d(\Z)$, then $Q$ is called
\textit{$G$-perfect} if for the cone $\Delta_F$ with $Q \in \Delta_F$
we have $\dim (\Delta_F \cap \MF(G)) = 1$. Two $G$-perfect forms $Q$
and $Q'$ are called \textit{$G$-equivalent} if there is a $g \in
N(\MB(G))$ and a positive $\lambda$ so that $g^t Q g = \lambda
Q'$. Otherwise they are called \textit{$G$-inequivalent}. By Voronoi's
first memoir \cite{voronoi-1907} the polyhedral subdivision $\MP$
satisfies the assertion of Theorem~\ref{th:equivariance}. Hence, there
are only finitely many $G$-inequivalent $G$-perfect forms by our 
Theorem \ref{th:equivariance}.

Usually $G$-perfect forms 
(cf. \cite{martinet-2003,bms-1992,jaquet-1995,opgenorth-1995,opgenorth-2001})
are defined via {\em normal cones} of faces $F_G$ of $\MP_m \cap \MF(G)$ in $\MF(G)$.
A face $F$ of $\MP_m$ is uniquely characterized by the
set $$\Min(F) = \{ \vec{v} \in \Z^d : Q[\vec{v}]=m \mbox{ for all } Q\in F \},$$
which is independent of the particular choice of $m$.
The normal cone $\MN(F)$ of $F$ is given by 
$\MN(F) = \cone \{ \vec{v}\vec{v}^t : \vec{v} \in \Min(F) \}$
and the normal cone of the face $F_G = F \cap \MF(G)$ in $\MF(G)$
is obtained by an orthogonal projection of $\MN(F)$ onto $\MF(G)$.
If different inner products are used, the resulting
cones differ as seen in the cases of \cite{jaquet-1995} and
\cite{opgenorth-1995}.
\end{remark}

\begin{remark}
\label{rem:closedgeodesics}
The case when the linear subspace $T$ is pointwise stabilized by a
finite subgroup of $\GL_d(\Z)$ is not the only case when there are
only finitely many $T$-inequivalent $T$-secondary
cones. Bayer--Fluckiger and Nebe (\cite[Theorem 3.1]{bn-2004}) give
another sufficient condition for $T$ which ensures
finiteness in the two-dimensional case. However, we are not aware of
further results regarding this question. 
For all the subspaces $T$ we considered for the applications, 
there were only finitely many $T$-inequivalent $T$-secondary cones. 
These subspaces are all spanned by rational forms, 
but this alone does not guarantee finiteness for $d\geq 3$,
as shown by the following example due to Yves Benoist (private communication):

Consider the subspace $T$ spanned by
$x^2+2y^2+z^2$ and $xy$. It contains the two 
positive semidefinite, but non-positive forms
$(x\pm\sqrt{2} y)^2 + z^2$.
In any small neighborhood of these two forms, we find
infinitely many forms in $T\cap \sdo$ with 
pairwise differing Delone subdivisions.
In particular, the minimal non-zero vectors $\vec{v}\in\Z^d$
tend towards the kernel of one of the two forms;
the Delone subdivisions have differing edges 
$\conv\{\vec{0},\vec{v}\}$ with $\vec{v}=(p,q,0)^t$
and $\tfrac{p}{q}$ sufficiently close to $\mp \sqrt{2}$.
Thus there are infinitely many different $T$-secondary cones,
but the setwise stabilizer of $T$ in $\gldz$ is finite. 
In order to see this, note that the three lines of
rank-$2$ forms in $T$ have to be permuted.
Note that the same argument does not apply to dimension $2$.
The space spanned by $x^2+2y^2$ and $xy$ has an infinite stabilizer.
\end{remark}


\section{Other generalizations of Voronoi's Reduction Theory}
\label{sec:comparison}

In the seventies, Mumford and Satake described a general procedure for
compactifying quotients of bounded symmetric domains by arithmetic
groups (cf. \cite{amrt-1975}). Their constructions use special
decompositions of the domains of interest into rational polyhedra. As
an application of this general method, Namikawa (\cite{namikawa-1976})
constructed a compactification of Siegel modular varieties using
Voronoi's reduction theory. Siegel modular varieties are moduli spaces
of Abelian varieties. Over the complex numbers they arise as quotient
of the Siegel upper halfspace by an arithmetic group.
Namikawa's construction opened the possibility to
understand the boundary of the moduli space with help of
degenerations. Later algebraic geometers refined this work and they
used also generalizations of Voronoi's reduction theory for this. In
this section we want to review these generalizations and compare them to
ours. For more information about the geometry of Siegel modular
varieties we refer to the survey \cite{hs-2002} of Hulek and Sankaran.

Oda and Seshadri study in \cite{os-1979} a generalization of Delone
subdivisions which they call Namikawa decompositions: Let $E$ be a
Euclidean space where an orthogonal decomposition into subspaces $E'$
and $E''$ is given: $E = E' \perp E''$, and let $\Lambda \subseteq E$
be a lattice so that $\Lambda \cap E'$ is also a lattice. For a vector
$\psi \in E''$ the Namikawa decomposition is
\[
D_{\psi} = \{\pi'(D(x)) : x \in E' + \psi\},
\]
where $\pi' : E \to E'$ is the orthogonal projection onto $E'$ and
\[
D(x) = \mbox{conv}\{\xi_1, \ldots, \xi_r\},
\]
with vertices $\xi_i \in \Lambda$ closest to $x$, is a Delone
polytope of $\Lambda$. Then, Oda and Seshadri explain in
\cite[Proposition 2.3]{os-1979} how the Namikawa decomposition changes
when one varies the vector $\psi$ along $E''$.  Our generalization of
Voronoi's reduction theory goes into a different direction: We fix a
lattice $\Lambda$ in a real vector space of finite dimension and a
basis of $\Lambda$ and study how the Delone subdivision of $\Lambda$
changes when one varies the inner product along a given subspace
(in the space of Gram matrices).

Alexeev and collaborators (cf. \cite{abh-2002,alexeev-2002,alexeev-2004}
and the survey \cite[Section~5]{alexeev-2006}) 
consider so-called semi-Delaunay decompositions: Let
$M \subseteq \R^d$ be a lattice and $\Gamma \subseteq M$ be a
sublattice of finite index. A $\Gamma$-periodic polyhedral subdivision
of $\R^d$ is called a semi-Delaunay decomposition when it is the
projection of the lower hull (see Section \ref{ssec:bistellar}) of the
lifted points $(m, h(m))$ with $m \in M$ and $h : M \to \R$ a function
of the form $h(m) = q(m) + r(m)$ where $q : \R^d \to \R$ is a positive
semidefinite quadratic form and $r : M / \Gamma \to \R$ an arbitrary
function on the cosets. Delone subdivisions of \textit{rational
standard periodic sets}, i.e.\ sets of the form
\eqref{eqn:standard_periodic_set} with rational $t'_i$, are
covered by this construction: Take $\Gamma = \Z^d$ and an overlattice $M$
containing all $t'_i$. Define $r(t'_i + \Gamma) = 0$ and $r(m
+ \Gamma)$ large enough for all $m \in M$ not in the rational standard
periodic set.

Nevertheless, we believe that our treatment has its merits since it is
very explicit and can be immediately used for computations. Moreover,
our generalization of Voronoi's reduction theory to positive definite
quadratic forms with prescribed automorphism group, which is essential
for our applications, is not covered by the reviewed constructions.


\section{Optimizing lattice sphere coverings}
\label{sec:coverings}

The lattice covering problem is a classical problem in ``Geometry of
Numbers''.  Roughly speaking, the problem is concerned with the
determination of the most economical way to cover~$\R^d$. We deal with
it in Section~\ref{sec:thinfields} and~\ref{sec:newcoverings} where we
apply our generalization of Voronoi's theory. In this section we give
the necessary definitions and describe briefly the methods used for
the applications. For further reading and more background we refer the
interested reader to \cite{sv-2005}.

\subsection{Definitions}

We assume that $L$ is a lattice of full rank $d$, that is, there
exists a matrix $B \in \gldr$ with $L = B\Z^d$.  The
\textit{determinant} $\det(L) = |\det(B)| > 0$ of $L$ is well defined
and does not depend on the chosen basis.  Let $B^d=\{\vec{x}\in\R^d :
\|\vec{x}\| \leq 1 \}$ denote the solid unit sphere in $\R^d$.  The
Minkowski sum $L + \alpha B^d = \{\vec{v} + \alpha\vec{x} : \vec{v}
\in L, \vec{x} \in B^d\}$, with $\alpha \in \R_{>0}$, is called a
\textit{lattice covering} if $\R^d = L + \alpha B^d$. The
\textit{covering radius} $\mu(L)$ of $L$ is given by
\[
\mu(L)
=
\min \{ \mu : \mbox{$L + \mu B^d$ is a lattice covering} \}.
\]
For $\alpha \in \R$ we have
\[
\mu(\alpha L) = |\alpha|\mu(L)
\]
Thus, the \textit{covering density} 
\[
\Theta(L)=\frac{\mu(L)^d}{\det(L)}\cdot\kappa_d
,
\]
with $\kappa_d = \vol B^d$, is invariant with respect to scaling
of~$L$.  Note also that $\Theta$ is an invariant of the isometry
classes.

The \textit{lattice covering problem} asks to minimize $\Theta$ among
all $d$-dimensional lattices.  It has been solved only for dimensions
$d\leq 5$ based on Voronoi's classical reduction theory and the
knowledge of all Delone subdivisions in these dimensions.

If we work with PQFs, the definitions above translate in the following way:
\[
\Theta(Q) = \Theta(L) = \sqrt{\frac{\mu(Q)^d}{\det Q}} \cdot \kappa_d,
\]
Here, the \textit{inhomogeneous minimum} $\mu(Q)$ is given by
\[
\mu(Q) = \max_{\vec{x} \in \R^d} \min_{\vec{v} \in \Z^d} Q[\vec{x} - \vec{v}],
\]
satisfying $\det(L) = \sqrt{\det(Q)}$, $\mu(L)=\sqrt{\mu(Q)}$, for a
corresponding lattice $L$ obtained from $Q$.

Given a linear subspace $T\subseteq\sd$, let  
\[
\Theta_T = \inf_{Q\in T\cap\sdo} \Theta (Q)
\]
denote the bound on the covering density with respect to $T$. This
bounds is computed in our applications described in Section
\ref{sec:thinfields} and Section \ref{sec:newcoverings}.

\subsection{Determinant maximization problems}

Let $T \subseteq \sd$ be a linear subspace.  To compute $\Theta_T$, we
consider all inequivalent $T$-generic Delone subdivision $\MD$ (if
possible) and minimize $\Theta$ among all PQFs in
$\overline{\SC_T(\MD)}$.  This can be achieved by solving a
\textit{determinant maximization problem}. These are convex
programming problems of the form
\begin{equation}   \label{eqn:primal-problem}
\begin{array}{ll}
\mbox{\textbf{minimize}} & \vec{c}^t\vec{x}-\log\det G(\vec{x})\\
\mbox{\textbf{subject to}} & G(\vec{x}) \succ 0, 
                             F(\vec{x}) \succeq 0,
\end{array}
\end{equation} 
over the variable vector $\vec{x} \in \R^D$.  The objective function
contains a linear part given by $\vec{c} \in \R^D$. The affine maps
$G : \R^D \to \R^{m \times m}$, as well as $F : \R^D \to \R^{n \times n}$, 
are given by
$$
\begin{array}{rcl}
G(\vec{x}) & = & G_0 + x_1 G_1 + \cdots + x_D G_D,\\
F(\vec{x}) & = & F_0 + x_1 F_1 + \cdots + x_D F_D,
\end{array}
$$ 
where $G_i \in \R^{m \times m}$ and $F_i \in \R^{n \times n}$, $i =
0,\ldots, D$, are symmetric matrices.  The notation $G(\vec{x}) \succ
0$ and $F(\vec{x}) \succeq 0$ gives the constraints ``$G(\vec{x})$ is
positive definite'' and ``$F(\vec{x})$ is positive semidefinite''.
Note that we specialize to a
\textit{semidefinite programming problem}, if $G(\vec{x})$
is the identity matrix for all $\vec{x} \in \R^D$.
One nice feature of determinant maximization problems is that there is
a duality theory similar to the one of linear programming (see \cite{vbw-1998}).
It allows us to compute an interval in which the optimum 
is attained, the so called \textit{duality gap}.

We can express the condition $\mu(Q)\leq 1$ as a \textit{linear matrix
inequality} (LMI) $F(\vec{x}) \succeq 0$ where the optimization vector
$\vec{x}$ is given by the coefficients of $Q = \sum_{i=1}^{\dim T}
x_i A_i$, with respect to some basis $(A_1, \dots, A_{\dim T})$
of $T$.  This is seen by the following proposition 
due to Delone et al \cite{ddrs-1970}
(cf. \cite[Proposition 7.1]{sv-2005}), together with the crucial
observation that the inner product $(\cdot, \cdot)$ defined by
$(\vec{y}, \vec{z}) = \vec{y}^tQ\vec{z}$ can be expressed as
\[
(\vec{y}, \vec{z}) = \langle Q , \vec{y} \vec{z}^t \rangle
 = \sum_{i=1}^{\dim T} x_i \langle A_i, \vec{y} \vec{z}^t \rangle
,
\]
hence as a linear combination of the parameters $x_i$.

\begin{proposition}
\label{prop:lmi}
Let $L = \conv\{\vec{0}, \vec{v}_1, \ldots, \vec{v}_d\}\subseteq
\R^d$ be a $d$-dimensional simplex. 
Then the circumradius of $L$ is at most~$1$ with respect to $(\cdot, \cdot)$ 
if and only if
$$
\BR_L(Q) =
\begin{pmatrix}
4 & (\vec{v}_1, \vec{v}_1) & (\vec{v}_2, \vec{v}_2) & \ldots &
(\vec{v}_d, \vec{v}_d)\\
(\vec{v}_1, \vec{v}_1) & (\vec{v}_1, \vec{v}_1) & (\vec{v}_1, \vec{v}_2) & \ldots &
(\vec{v}_1, \vec{v}_d)\\
(\vec{v}_2, \vec{v}_2) & (\vec{v}_2, \vec{v}_1) & (\vec{v}_2,
\vec{v}_2) & \ldots & (\vec{v}_2, \vec{v}_d)\\
\vdots & \vdots & \vdots & \ddots & \vdots\\
(\vec{v}_d, \vec{v}_d) & (\vec{v}_d, \vec{v}_1) & (\vec{v}_d, \vec{v}_2) & \ldots &
(\vec{v}_d, \vec{v}_d)\\
\end{pmatrix}
\succeq 0.
$$
\end{proposition}

Since a block matrix is semidefinite if and only if the blocks are
semidefinite, we get the desired result. Here, we call two Delone polytopes
$P$ and $P'$ \textit{equivalent with respect to $Q \in \sdo$}, if there exists a 
$U$ in the \textit{automorphism group} $$\Aut(Q)=\{U\in\gldz : U^t Q U = Q \}$$ of $Q$ and a $\vec{v}\in\Z^d$ with 
$P=\vec{v}+UP'$.

\begin{proposition}  
\label{propos:mu}
Let $Q \in \sdo$ be a PQF with Delone subdivision $\MD$. Let $P_1,
\ldots, P_n$ be a representative system of $d$-dimensional 
Delone polytopes in $\MD$, which are inequivalent with respect to $Q$. 
For every $P_i$ choose a $d$-dimensional
simplex $L_i$ with $\vertex L_i \subseteq \vertex P_i$. Then
\[
\mu(Q)\leq 1
\quad\Longleftrightarrow\quad
\begin{pmatrix}
\boxed{\BR_{L_1}(Q)} & 0                   & 0 & \ldots & 0\\
0                   & \boxed{\BR_{L_2}(Q)} & 0 & \ldots & 0\\
\hdotsfor{5}\\
0                   & 0                   & 0 & \ldots & \boxed{\BR_{L_n}(Q)}
\end{pmatrix}
\succeq 0.
\]
\end{proposition}

Thus the constraint $\mu(Q)\leq 1$ can be brought into one LMI of type
$F(\vec{x}) \succeq 0$. We can add linear constraints on the
parameters $x_i$ by extending $F$ by a $1\times 1$ block matrix for
each linear inequality. In this way we obtain one LMI for the two
constraints $\mu(Q)\leq 1$ and $Q \in
\overline{\SC(\MD)}$. 
Therefore we can determine a PQF in $\overline{\SC_T(\MD)}$ minimizing
the covering density by solving the following determinant maximization
problem.
$$
\begin{array}{ll}
\mbox{\textbf{minimize}} & -\log\det(Q)\\
\mbox{\textbf{subject to}} & \mbox{$Q \in \overline{\SC_T(\MD)}$, $\mu(Q) \leq 1$}.\\
\end{array}
$$ 
Here $F(\vec{x})$ is as described above, and $\vec{c} = \vec{0}$, and
$G(\vec{x}) = \sum_{i = 1}^{\dim T} x_i A_i$ in
(\ref{eqn:primal-problem}) where the $A_i$ form a basis of $T$.

\begin{remark}
The problem is a convex programming problem. Hence (see e.g.\
\cite[Proposition 9.1]{sv-2005}), if the automorphism group of the
considered Delone subdivision $\MD$ contains the pointwise stabilizer
of the subspace $T$, we know that a PQF minimizing the covering
density in $\overline{\SC(\MD)}$ is contained in $T$ . Hence by
choosing the subspace $T$ carefully we sometimes can dramatically
decrease the dimension of the optimization problem.
\end{remark}

\begin{remark}
Determinant maximization problems have not only theoretically nice
features, but also can be solved in practice quite well, for example
with interior point methods.  For the solution of the optimization
problems described above we use the software package \texttt{MAXDET}
\cite{maxdet-1996} of Wu, Vandenberghe, and Boyd as a subroutine.  By
using rational approximations it is possible to give mathematical
rigorous error bounds or even proofs of local optimality.  For details
we refer to \cite{sv-2005}.  Utilizing these ideas we developed a
program \texttt{rmd} (\textit{rigorous maxdet}) and based on it
\texttt{coop}, a \textit{covering optimizer}.  These programs,
together with a short tutorial can be obtained from our web page
\cite{latgeo-2005}.  They allow to approximate --- mathematically
rigorously --- the optimum covering density with respect to a given
$T$-secondary cone.
\end{remark}

\section{Algorithmic issues}
\label{sec:algorithms}

In this section we discuss several algorithmic issues which have to be
resolved to turn our generalization of Voronoi's reduction theory into
an effective procedure. As in Section \ref{sec:equivariant} we assume that
all Delone subdivisions have $\Z^d$ as a vertex-set.
For the original theory of Voronoi a similar
discussion can be found in \cite[Section 5.3]{sv-2005}. Here we
emphasize those points where the generalization differs from the
classical case. First we adapt \cite[Algorithm 1]{sv-2005}, which
enumerates all inequivalent Delone triangulations, to an algorithm
enumerating all $T$-inequivalent $T$-generic Delone subdivisions.
As in Example \ref{ex:infinite-many}, our Algorithm 1 below does not necessarily 
stop after finitely many steps, depending on $T$.

Apart from this, when we compare Algorithm 1 to the one
of the original theory, essentially two new tasks arise: 
\begin{enumerate}
\item We have to find an initial $T$-generic Delone subdivision
(cf. Algorithm 2), whereas in the classical case this can be given by
``Voronoi's first subdivision''. 
\item We have to check whether two
$T$-generic Delone subdivisions are $T$-equivalent (cf. Algorithm 3).
\end{enumerate}

\begin{figure}[htb]
\begin{center}
\fbox{
\begin{minipage}{12.5cm}
\begin{flushleft}
\smallskip
\textbf{Input:} Linear subspace $T \subseteq \MS^d$\\
\textbf{Output:} Set~$\MNR$ of all $T$-inequivalent $T$-generic Delone subdivisions\\
\smallskip
$Y \leftarrow \{\MD_1\}$, where $\MD_1$ is a $T$-generic Delone subdivision (Algorithm 2).\\
$\MNR \leftarrow \emptyset$.\\
\textbf{while} there is a $\MD \in Y$ \textbf{do}\\
\hspace{2ex} $Y \leftarrow Y \backslash \{\MD\}$. $\MNR \leftarrow \MNR \cup \{\MD\}$.\\
\hspace{2ex} Compute the linear inequalities of $\overline{\SC_T(\MD)}$ as in Theorem~\ref{th:secondarycone}.\\
\hspace{2ex} Compute the facets $F_1, \ldots, F_n$ of $\overline{\SC_T(\MD)}$.\\
\hspace{2ex} \textbf{for} $i = 1, \ldots, n$ \textbf{do}\\
\hspace{2ex}\hspace{2ex} \textbf{if} $F_i$ is not a dead-end \textbf{then}\\
\hspace{2ex}\hspace{2ex}\hspace{2ex} Compute the bistellar neighbor $\MD_i$ of $\MD$, defined by $F_i$ as in Theorem~\ref{th:bistellar}.\\
\hspace{2ex}\hspace{2ex}\hspace{2ex}  \textbf{if} $\MD_i$ is not $T$-equivalent to a Delone subdivision in the set\\
\hspace{2ex}\hspace{2ex}\hspace{2ex} \phantom{\textbf{if}}
$\MNR \cup \{\MD_j : \mbox{$j \in \{1, \ldots, i-1\}$, $F_j$ not a dead-end}\}$ (Algorithm 3) \textbf{then}\\
\hspace{2ex}\hspace{2ex}\hspace{2ex}\hspace{2ex} $Y \leftarrow Y \cup \{\MD_i\}$.\\
\hspace{2ex}\hspace{2ex}\hspace{2ex} \textbf{end if}\\
\hspace{2ex}\hspace{2ex} \textbf{end if}\\
\hspace{2ex} \textbf{end for}\\
\textbf{end while}\\
\end{flushleft}
\end{minipage}
}
\\[1ex]
\textbf{Algorithm 1.} Enumeration of all $T$-inequivalent $T$-generic Delone subdivisions.
\end{center}
\end{figure}

For Algorithm 2 we need one more definition: The dimension of the
linear span of $\SC_T(\MD)$ is called \textit{$T$-rigidity index}. For
the classical setting $T = \MS^d$ the \textit{rigidity index} was
introduced by Baranovskii and Grishukhin in \cite{bg-2002}. Algorithm
2 is a randomized algorithm and belongs to the class of so-called Las
Vegas algorithms. That is, it always produces correct results whereas
the running time is a random variable.

\begin{figure}[htb]
\begin{center}
\fbox{
\begin{minipage}{12.5cm}
\begin{flushleft}
\smallskip
\textbf{Input:} Linear subspace $T \subseteq \MS^d$\\
\textbf{Output:} A $T$-generic Delone subdivision $\MD_1$\\
\smallskip
\textbf{repeat}\\
\hspace{2ex} Choose a random $Q$ in $T \cap \sdo$.\\
\hspace{2ex} Compute $\MD_1 = \Del(Q)$ (see \cite{dsv-2006}).\\
\hspace{2ex} Compute the $T$-rigidity index $m$ of $\MD_1$ by the linear equalities in Theorem~\ref{th:secondarycone}.\\
\textbf{until} $m = \dim T$.\\
\end{flushleft}
\end{minipage}
}
\\[1ex]
\textbf{Algorithm 2.} Finding a $T$-generic Delone subdivision.
\end{center}
\end{figure}

Algorithm 3 checks whether two $T$-generic $T$-secondary cones are
$T$-equivalent. For this we need a definition.

\begin{definition}
We say that a positive
semidefinite quadratic form $Q$ is \textit{rational normalized} if it
is integral and the entries have greatest common divisor equal to one.
\end{definition}

For Algorithm 3 we have to test whether two PQFs are equivalent and we
have to compute the automorphism group of a given PQF. For both tasks
exist efficient algorithms by Plesken and Souvignier
\cite{ps-1997}. Implementations are part of the computational algebra
system \texttt{MAGMA} and of the computer package \texttt{CARAT}
(Crystallographic AlgoRithms And Tables, cf.~\cite{ops-1998}).

\begin{figure}[htb]
\begin{center}
\fbox{
\begin{minipage}{12.5cm}
\begin{flushleft}
\smallskip
\textbf{Input:} Linear subspace $T \subseteq \MS^d$, and $T$-generic Delone subdivisions $\MD_1$, $\MD_2$\\
\textbf{Output:} Yes, if $\MD_1$ and $\MD_2$ are $T$-equivalent; No, otherwise\\
\smallskip
$\MNR_i \leftarrow \{\mbox{rational normalized element $R$ of extreme rays of
$\overline{\SC_T(\MD_i)}$}\}$, for $i = 1, 2$.\\ 
$Q_i \leftarrow \sum_{R \in \MNR_i} R$, for $i = 1, 2$. \\
Compute $\Aut(Q_1)$.\\
\textbf{if} there exists $A \in \GL_d(\Z)$ and $B \in \Aut(Q_1)$\\ 
\hspace{2ex} with $A^t Q_1 A = Q_2$ and $(BA)^t T (BA) = T$ \textbf{then}\\
\hspace{2ex} \textbf{return} Yes. \\
\textbf{else}\\
\hspace{2ex} \textbf{return} No.\\
\textbf{end if}\\
\end{flushleft}
\end{minipage}
}
\\[1ex]
\textbf{Algorithm 3.} Checking $T$-equivalence.
\end{center}
\end{figure}

The closed polyhedral cones $\overline{\SC_T(\MD)}$ are
\textit{rational} (assuming that the vertex-set is $\Z^d$).  That is,
they have a description by rational inequalities (cf. Section
\ref{sec:generalization}) as well as a description
$\overline{\SC_T(\MD)} = \cone \{ R_1,\dots, R_k\}$ as a cone
generated by finitely many rational normalized quadratic forms lying
in extreme rays of $\overline{\SC_T(\MD)}$.

Dealing with symmetries of a polyhedral cone, the
\textit{characteristic form} $Q=\sum_{i = 1}^k R_i$ is of great
importance.  This is due to the following proposition which is used in
Algorithm 3.  We call two rational polyhedral cones $\MP$ and $\MP'$
in $\sd$ \textit{equivalent}, if there exists a $U\in\gldz$ such that
$U^t \MP U = \MP'$. As in the proof of Theorem~\ref{th:equivariance},
the {\em automorphism group} of $\MP$ is defined by $\Aut(\MP) = \{ U
\in \gldz : U^t \MP U = \MP \}$.

\begin{lemma}
The automorphism group of a characteristic form contains the
automorphism group of its polyhedral cone.  The characteristic forms
of two equivalent rational polyhedral cones are equivalent.
\end{lemma}

\begin{proof}
Suppose that for two polyhedral cones $\MP,\MP'\subset \sdclosed$
there exists a $U\in\gldz$ such that $U^t \MP U = \MP'$.
Extreme rays of $\MP$ are mapped onto extreme rays of $\MP'$.  We need
to show that this is also true for the uniquely determined rational
normalized elements of the extreme rays. For this let $Q$ be such a rational
normalized element of an extreme ray of $\MP$. It is mapped onto the integral
matrix $U^tQU=\alpha Q'$ defining an extreme ray of $\MP'$, with rational
normalized $Q'$. Thus $\alpha\in\N$.  On the other hand $Q'$ is mapped
onto an integral multiple of $Q$ via $U^{-1}$ and $(U^{-1})^t
Q'(U^{-1}) = (1/\alpha) Q$. Thus $\alpha=1$, which proves the
assertion.
\end{proof}

\section{Application I: Classification of totally real thin number fields}
\label{sec:thinfields}

In this section we apply our theory to a problem in algebraic number
theory. Recently, Bayer--Fluckiger introduced in \cite{bayer-2004} the
notion of thin algebraic number fields. A thin algebraic number field
is Euclidean due to a special geometric reason, which we explain in
Section~\ref{ssec:thinbackground}. She proved that there exist only
finitely many thin number fields (\cite[Proposition~11.4]{bayer-2004}).

In \cite{bn-2004} Bayer--Fluckiger and Nebe gave a complete list
(\cite[Theorem 5.1]{bn-2004}) of 17 candidates for totally real thin
algebraic number fields. They proved that 13 of them are thin and one
of them is weakly thin but not thin. However, in three cases they did
not know how to decide whether the fields are thin or not. Using our
theory we finish the classification of totally real thin algebraic
number fields. In particular we show in Section~\ref{ssec:thinresults}
that the three open cases do not give weakly thin fields.

\subsection{Background and definitions}
\label{ssec:thinbackground}

Let us recall some standard definitions from algebraic number theory.

Let $K$ be a number field of degree $n = [K : \Q]$. From now on we
 assume that $K$ is \textit{totally real}, that is, for all embeddings
 of fields $\sigma_i : K \to \C$ with $i = 1, \ldots, n$ we have
 $\sigma_i(K) \subseteq \R$.

Let $\oK$ be its ring of integers which is a $\Z$-module of rank
$n$. Thus we can write $\oK = \Z\omega_1 + \cdots + \Z\omega_n$. We
embed $\oK$ into $\R^n$ via the map $\sigma$ defined by
$\sigma(x) = (\sigma_1(x), \ldots, \sigma_n(x))$. Hence, $\sigma(\oK)$
is a lattice of rank $n$.

By $\MP$ we denote the set of $\alpha \in K$ with $\sigma_i(\alpha) >
0$ for all $i = 1, \ldots, n$. Then
\[
\langle x, y \rangle_{\alpha} = \sum_{i = 1}^n \sigma_i(\alpha x y)
= \sum_{i = 1}^n \sigma_i(\alpha) \sigma_i(x) \sigma_i(y)
\]
defines an inner product on $K$ and so on $\R^n$. Every inner product
$\langle \cdot, \cdot \rangle : K \times K \to \R$ with the additional
property $\langle x, yz \rangle = \langle xz, y\rangle$
for all $x,y,z \in K$ is of this form.

We denote by $\Lambda_{K,\alpha}$ the lattice which is given by the
pair $(\oK, \langle \cdot, \cdot \rangle_{\alpha})$ where $\alpha \in
\MP$.  We say that $K$ is \textit{weakly thin} if for $\Theta(K) =
\min_{\alpha \in \MP} \Theta(\Lambda_{K,\alpha})$ we have
the inequality
\[
\Theta(K) \leq \sqrt{\frac{n^n}{(\det \Lambda_{K,1})^2}} \cdot \kappa_n
,
\]
and we say that $K$ is \textit{thin} if we have strict inequality. By
$t(K)$ we denote the number on the right hand side of this inequality.

Let us briefly explain the motivation of this definition. First note
that $(\det \Lambda_{K,1})^2$ equals the discriminant $d_K$ of
$K$. Let $N : K \to \R$ be the norm of $K$ which is given by $N(x) =
\prod_{i=1}^n \sigma_i(x)$. A number field $K$ is called
\textit{Euclidean} if its ring of integers $\oK$ is an Euclidean ring
with respect to the absolute value of the norm function. Equivalently,
$K$ is Euclidean if and only if its Euclidean minimum
\[
M(K) = \sup_{x \in K} \inf_{y \in \oK} |N(x-y)|
\]
is strictly less than $1$. Currently there is no algorithm known which
computes the Euclidean minimum of a number field. It is also an open
problem if there exists a finite or an infinite number of Euclidean
number fields. In \cite{lenstra-1976} Lenstra, motivated by work of
Hurwitz, gave several bounds on the Euclidean minimum of a number
field using methods from ``Geometry of Numbers''. In \cite{bayer-2004}
Bayer-Fluckiger extended these results. By using the inequality
between the arithmetic and geometric means to relate the norm function
$N$ to the norm $\|x\|_{\alpha} = \sqrt{\langle x, x \rangle_{\alpha}}
= \sqrt{\sum_{i=1}^n \sigma_i(\alpha x^2)}$, she showed that thin
fields are Euclidean fields (\cite[Proposition 11.2]{bayer-2004}) and
that there exist only finitely many thin fields (\cite[Proposition 11.4]{bayer-2004}).

\subsection{Classification}
\label{ssec:thinresults}

By applying the lower bound for sphere coverings of Coxeter, Few and
Rogers \cite{cfr-1959}, Bayer-Fluckiger showed that the degree of a
totally real thin number field is at most $5$ (\cite[Proof of
Proposition 11.4]{bayer-2004}). All number fields having low degree
and low discriminant are known (\cite{kant-2005}). Using this list
Bayer--Fluckiger and Nebe gave a complete list of $17$ candidates of
totally real number fields which might be thin. In Table~1 we list
these candidates $K$ together with the relevant parameters: The degree
$n$, the discriminant $d_K$, the bound $t(K)$ we defined above and the
minimum covering density $\Theta(K)$.

\begin{figure}
\begin{center}
\begin{tabular}{c|c|c|c|c|c}
$n$ & $d_K$ & $K$ & $t(K)$ & $\Theta(K)$ & thin\\
\hline
$2$ & $5$ & $\Q[x]/(x^2 - 5)$   & $\geq 2.8099$  & $\leq 1.2645$ & yes \\
    & $8$ & $\Q[x]/(x^2 - 2)$   & $\geq 2.2214$  & $\leq 1.3463$ & yes \\
    & $12$ & $\Q[x]/(x^2 - 3)$  & $\geq 1.8137$  & $\leq 1.2092$ & yes \\
    & $13$ & $\Q[x]/(x^2 - 13)$ & $\geq 1.7426$  & $\leq 1.5708$ & yes \\
    & $17$ & $\Q[x]/(x^2 - 17)$ & $\geq 1.5238$  & $\leq 1.2497$ & yes \\
    & $21$ & $\Q[x]/(x^2 - 21)$ & $\geq 1.3711$  & $\leq 1.2242$ & yes \\
    & $24$ & $\Q[x]/(x^2 - 6)$  & $\geq 1.2825$  & $\leq 1.2583$ & yes \\
\hline
$3$ & $49$ & $\Q[x]/(x^3 + x^2 -2x - 1)$   & $\geq 3.1093$ & $\leq 1.5584$ & yes \\
    & $81$ & $\Q[x]/(x^3 - 3x + 1)$        & $\geq 2.4183$ & $\leq 2.1225$ & yes \\
    & $148$ & $\Q[x]/(x^3 + x^2 - 3x - 1)$ & $\geq 1.7891$ & $\leq 1.7014$ & yes \\
    & $169$ & $\Q[x]/(x^3 + x^2 - 4x + 1)$ & $\leq 1.6743$ & $\geq 1.8544$ & no \\
\hline
$4$ & $725$ & $\Q[x]/(x^4 - x^3 - 3x^2 + x + 1)$    & $\geq 2.9323$ & $\leq 2.7045$ & yes \\
    & $1125$ & $\Q[x]/(x^4 - x^3 - 4x^2 + 4x + 1)$  & $\geq 2.3540$ & $\leq 2.2935$ & yes \\
    & $1600$ & $\Q[x]/(x^4 + 2x^3 - 5x^2 - 6x - 1)$ & $\leq 1.9740$ & $\geq 2.2853$ & no \\
    & $1957$ & $\Q[x]/(x^4 - 4x^2 - x + 1)$         & $\leq 1.7849$ & $\geq 1.8939$ & no \\
    & $2000$ & $\Q[x]/(x^4 - 5x^2 + 5)$             & $= \Theta_4$ & $= \Theta_4$ & weakly \\
\hline
$5$ & $14641$ & $\Q[x]/(x^5 + x^4 - 4x^3 - 3x^2 + 3x + 1)$ & $\geq 2.4318$ & $\leq 2.2961$ & yes
\end{tabular}
\par\smallskip
\textbf{\textsf{Table 1.}} Classification of totally real thin number fields.
\end{center}
\end{figure}

Bayer--Fluckiger and Nebe showed by giving a number $\alpha$, which
defines the appropriate inner product $\langle \cdot, \cdot
\rangle_{\alpha}$, that $13$ candidates are thin. 
For computing an upper bound of $\Theta(K)$ as given in the table
we used their values $\alpha\in\MP$ in all but one case.
In the case of $\Q[x]/(x^3 + x^2 - 3x - 1)$ their table contains a
misprint. Instead of $\alpha = 1 - 18\bar{x} + 10\bar{x}^2$ which does not lie
in $\MP$ we use $\alpha = 2-\bar{x}$ instead.
For one candidate Bayer--Fluckiger and Nebe showed that it is weakly thin. 
In this case $t(K) = \Theta(K) = \Theta_4$ where $\Theta_4$ is the least 
lattice covering density among all $4$-dimensional lattices. 
It is uniquely attained by the lattice $\mathsf{A}^*_4$.

The other three cases were left open by Bayer--Fluckiger and Nebe 
and by our computation it turns
out that they do not give thin fields. Generally, for a totally real
number field $K$, there exists a subspace $T$ of $\MS^n$ 
whose dimension equals the degree of $K$, and such that $\Theta(K)= \Theta_T$.
The corresponding subspace $T$ 
is given by the basis $\begin{pmatrix} \langle \omega_i,
\omega_j \rangle_{\alpha_k}
\end{pmatrix}_{1 \leq i, j \leq n}
$, $k = 1, \ldots, n$,
where $\oK = \Z \omega_1 + \cdots + \Z \omega_n$ and $\alpha_1,
\ldots, \alpha_n$ forms a $\Q$-basis of $K$ with $\alpha_k \in \MP$.
Thus if there exist only finitely many inequivalent $T$-secondary cones,
we can compute $\Theta(K) = \Theta_T$ by the
methods explained in Section~\ref{sec:coverings}. 

For the field $K = \Q[x]/(x^3 + x^2 - 4x + 1)$ we have $7$ generic
$T$-secondary cones and for all lattices $\Lambda_{K,\alpha}$ our
computation proves the bound $\Theta(\Lambda_{K,\alpha}) \geq
1.8544$. For the field $K = \Q[x]/(x^4 + 2x^3 - 5x^2 - 6x - 1)$ we
have $47$ generic $T$-secondary cones and for all lattices
$\Lambda_{K,\alpha}$ our computation shows $\Theta(\Lambda_{K,\alpha})
\geq 2.2853$. For the field $K = \Q[x]/(x^4 - 4x^2 - x + 1)$ we have
$341$ generic $T$-secondary cones and for all lattices
$\Lambda_{K,\alpha}$ our computation gives $\Theta(\Lambda_{K,\alpha})
\geq 1.8939$.

The fact that in all these cases the number of inequivalent
$T$-generic secondary cones is finite comes as a pleasant surprise.
In this situation $T$ is not the space of invariant forms of a Bravais
group so that we can not a priori rely on the finiteness result of
Section~\ref{sec:equivariant}. It remains an interesting open problem
to give a proof of this fact, which does not rely on the complete enumeration
of $T$-secondary cones.


\section{Application II: New best known sphere coverings}
\label{sec:newcoverings}

In this section we explain how the new theory can be used to construct
good and even new best known sphere coverings. With the described
methods, we were in particular able to construct new best known
(lattice) sphere coverings in dimensions $9,\dots, 15$. Table 2 gives
an overview on the currently best known sphere coverings up to
dimension $24$.  This table is an update of the table given in
\cite[Table 2.1]{cs-1988}.  Note in particular, that in comparison to the table
there, we have new best known sphere coverings in all dimensions $d
\in \{6,\dots,21\}\setminus\{16,18\}$.  Note that the problem has been
solved for $d\leq 5$ only, by the work of Ryshkov and Baranovskii \cite{rb-1975}.  
We keep an updated list with additional
informations on the involved lattices on our web page
\cite{latgeo-2005}.

This section is organized as follows: First we give some background
information on the lattices of Table 2.  Then the next two sections
deal with two different constructions of subspaces $T$ we used. Both
contain specific information on how we obtained the new best known
lattices. The last section gives some information on a similar
approach to the closely related packing-covering problem.

\begin{figure}
\begin{center}
\begin{tabular}{c|c|c||c|c|c}
$d$ & lattice & covering density $\Theta$ & $d$ & lattice & covering density $\Theta$ \\  
\hline 
 $1$ & $\Z^1$               &  $1$          & $13$ & ${\mathsf L}^{c}_{13}$ & $7.762108$ \\ 
 $2$ & ${\mathsf A}^*_2$    &  $1.209199$   & $14$ & ${\mathsf L}^c_{14}$ & $8.825210$ \\ 
 $3$ & ${\mathsf A}^*_3$    &  $1.463505$   & $15$ & ${\mathsf L}^c_{15}$ & $11.004951$ \\
 $4$ & ${\mathsf A}^*_4$    &  $1.765529$   & $16$ & ${\mathsf A}^*_{16}$ & $15.310927$ \\
 $5$ & ${\mathsf A}^*_5$    &  $2.124286$   & $17$ & ${\mathsf A}^9_{17}$ & $12.357468$ \\ 
 $6$ & ${\mathsf L}^{c}_6$  &  $2.464801$   & $18$ & ${\mathsf A}^*_{18}$ & $21.840949$ \\ 
 $7$ & ${\mathsf L}^{c}_7$  &  $2.900024$   & $19$ & ${\mathsf A}^{10}_{19}$ & $21.229200$ \\
 $8$ & ${\mathsf L}^{c}_8$  &  $3.142202$   & $20$ & ${\mathsf A}^7_{20}$ & $20.366828$ \\ 
 $9$ & ${\mathsf L}^{c}_9$  &  $4.268575$   & $21$ & ${\mathsf A}^{11}_{21}$ & $27.773140$ \\ 
$10$ & ${\mathsf L}^{c}_{10}$ & $5.154463$  & $22$ & $\Lambda_{22}^*$     & $\leq 27.8839$ \\
$11$ & ${\mathsf L}^{c}_{11}$ & $5.505591$  & $23$ & $\Lambda_{23}^*$     & $\leq 15.3218$ \\
$12$ & ${\mathsf L}^{c}_{12}$ & $7.465518$  & $24$ & $\Lambda_{24}$       & $7.903536$ \\ 
\end{tabular}
\par\smallskip
\textbf{\textsf{Table 2.}} Least dense known (lattice) coverings
up to dimension $24$.
\end{center}
\end{figure}

\subsection{Some background on best known lattice coverings}

The new sphere coverings in dimension $6$, $7$ and $8$ 
were obtained and described in detail in \cite{sv-2005} and \cite{sv-2005b}.
The lattices ${\mathsf A}^r_d$, where $r$ divides $d+1$, are the {\em Coxeter lattices} \cite{coxeter-1951}. 
A possible definition 
is via the root lattice
$$
{\mathsf A}_d= \{ \vec{x} \in \Z^{d+1} : \sum_{i=0}^{d+1} x_i=0 \}
.
$$
The Coxeter lattice is the lattice generated by ${\mathsf A}_d$
and the vector
$$
(1/r) \left( \sum_{i=1}^{d+1} \vec{e}_i \right)  - (\vec{e}_1+....+\vec{e}_{(d+1)/r} ) 
,
$$
where $\vec{e}_i$ denotes the $i$-th standard basis vector.
It is the unique sublattice of ${\mathsf A}_d^\ast$
containing ${\mathsf A}_d$ as a sublattice of index $r$.
In particular ${\mathsf A}_d^\ast={\mathsf A}_d^{d+1}$.
Other well known lattices in the series are 
${\mathsf A}_7^2 = {\mathsf E}_7$, ${\mathsf A}_7^4={\mathsf E}_7^*$ 
and ${\mathsf A}_8^3={\mathsf E}_8$.

Since the symmetric group $\mathsf{S}_{d+1}$ acts on the lattice
${\mathsf A}^r_d$, it is possible (with the help of a computer) to
enumerate all its orbits of Delone polytopes and to compute their
covering densities in fairly large dimensions. A detailed description
of the method together with results up to dimension $27$ can be found
in \cite{dsv-2006}.

Baranovskii \cite{baranovskii-1994} computed (by hand) the Delone
decomposition of ${\mathsf A}^5_9$, finding the former best known sphere
covering in dimension $9$.  Anzin \cite{anzin-2002} computed the
Delone decomposition of ${\mathsf A}_{11}^4$ and ${\mathsf A}_{13}^7$
establishing the former covering records in those dimensions.  In a
private communication he reported on computing the covering densities
of ${\mathsf A}_{14}^5$ and ${\mathsf A}_{15}^8$, which were also
best known ones in their dimension. Hence suitable Coxeter lattices provide
good covering lattices.  Nevertheless, by applying our theory, we
found that all of the five mentioned lattices do not even give a
locally optimal lattice covering.  However, these lattices give good
starting points. In fact, we obtained the new covering records in
dimensions $d=9,11,13,14,15$ by applying our theory to a suitable
linear subspace $T$ containing a PQF of the corresponding Coxeter lattices. 
More details are given in the next section.  For dimensions $d=10,12$ we 
used a lamination technique, which is described thereafter.

The entries of Table 2 for dimensions $d=22,23,24$ are ``consequences'' 
of the existence of the Leech lattice $\Lambda_{24}$.  
It may not surprise that the Leech lattice itself 
yields the best known lattice covering in dimension $24$.  
Its covering density was computed by Conway, Parker and Sloane
(\cite[Chapter 23]{cs-1988}). In \cite{sv-2005b} it is shown that
the Leech lattice gives at least a local optimum of the covering 
function $\Theta(Q)$. Note that the root lattice $\mathsf{E}_8$ 
does not have this property.
Knowing the comparatively low covering density of the Leech lattice,
Smith \cite{smith-1988} was able 
to estimate the covering densities of the {\em dual laminated lattices}
$\Lambda^*_{22}$ and $\Lambda^*_{23}$.

\subsection{Large subgroups and small linear subspaces}

One strategy to find good covering lattices is to consider subspaces $T$
containing a PQF, which gives a good or even best known sphere covering. 
In order to keep the number of $T$-secondary cones low (manageable),
one preferably chooses a subspace $T$ of low dimension, e.g. less than or equal to $4$.
Another problem is the size of Delone subdivisions
to be dealt with. If $T$ is contained in a space of invariant forms $\MF(G)$
of a preferably large finite subgroup $G$ of $\gldz$, then this problem
can be reduced by exploiting the symmetries of the Delone subdivisions, respectively
those of the forms $Q\in \MF(G)$ (see Proposition \ref{propos:mu}).

{\bf Dimension $\mathbf{9}$}. 
We optimized over a $3$-dimensional subspace $T$ containing a PQF of the former record lattice ${\mathsf A}_9^5$. 
We computed $210$ $T$-generic $T$-secondary cones. 
The new covering record is attained by a PQF which has $34$ orbits of Delone polytopes.                              

{\bf Dimension $\mathbf{11}$}. 
We optimized over a $3$-dimensional subspace $T$ containing a PQF of the former record lattice ${\mathsf A}_{11}^4$. 
We computed $2444$ $T$-generic $T$-secondary cones. 
The new covering record is attained by a PQF which has $99$ orbits of Delone polytopes.       

{\bf Dimension $\mathbf{13}$}. 
We optimized over a $2$-dimensional subspace $T$ containing a PQF of the former record lattice ${\mathsf A}_{13}^7$. 
We computed $79$ $T$-generic $T$-secondary cones. 
The new covering record is attained by a PQF which has $134$ orbits of Delone polytopes.

{\bf Dimension $\mathbf{14}$}.
We optimized over a $2$-dimensional subspace $T$ containing a PQF of the former record lattice ${\mathsf A}_{14}^5$.
We computed $162$ $T$-generic $T$-secondary cones. 
The new covering record is attained by a PQF which has $983$ orbits of Delone polytopes.

{\bf Dimension $\mathbf{15}$}. 
We optimized over a $2$-dimensional subspace $T$ containing a PQF of the former record lattice ${\mathsf A}_{15}^8$. 
We computed $109$ $T$-generic $T$-secondary cones. 
The new covering record is attained by a PQF which has $203$ orbits of Delone polytopes.

\subsection{Best coverings from laminations}

Another fruitful strategy is the construction of good coverings from
lower dimensional ones.  Assume that $Q\in\sdo$ is contained in a
linear subspace $T$ of $\sd$. We consider a corresponding lattice $L =
A\Z^d$ with $A\in\gldr$.  Further, we choose a point $A\vec{c}$ with
$\vec{c}\in\R^d$ and consider all lattices $L_{\lambda}$ in $\R^{d+1}$
generated by
$$
A'=
\begin{pmatrix}
A & A\vec{c}\\
0 & \lambda\\ 
\end{pmatrix}
\in\mathsf{GL}_{d+1}(\R)
$$
with $\lambda > 0$ and 
associated PQF
$$
Q'=\begin{pmatrix} Q & Q\vec{c} \\ \vec{c}^t Q & \vec{c}^t Q\vec{c}+\lambda^2 \end{pmatrix}
\in\MS_{>0}^{d+1}
.
$$ Choosing $Q$ within $T$ and $\lambda^2\in\R$, we obtain a linear
subspace $T'$ of forms $Q'$ with $\dim T' = \dim T + 1$. If $A\vec{c}$
is the center of a Delone polytope $P$ of the lattice $L$, then the
lattice vectors in the lattice $L_\lambda$ belonging to the $i$-th
layer
\[
\left\{
\begin{pmatrix}
A & A\vec{c}\\
0 & \lambda\\ 
\end{pmatrix}
\begin{pmatrix}
\vec{v}\\
i
\end{pmatrix}
: \vec{v} \in \Z^d
\right\}
,
\quad
\mbox{with $i \in \Z$,}
\]
have at least distance $|i|\lambda$ from $(A\vec{c}, 0)$. Hence, $P$
is a Delone polytope embedded in $L_{\lambda}$ when $\lambda$ is at
least the circumradius of $P$ in $L$.

{\bf Dimension $\mathbf{10}$}.  We considered the $1$-dimensional
linear space $T$ containing ${\mathsf A}_9^5$.  We look at the linear
subspace $T'$ of $\MS^{10}$ of dimension $2$ constructed with the
center of a Delone polytope having the largest symmetry group in
${\mathsf A}_9^5$.  We computed $6$ $T'$-generic $T'$-secondary cones.
The new covering record is attained by a PQF which has $4$ orbits of
Delone polytopes.

{\bf Dimension $\mathbf{12}$}.  We considered the $1$-dimensional
linear space $T$ containing ${\mathsf A}_{11}^4$.  We look at the
linear subspace $T'$ of $\MS^{12}$ of dimension $2$ constructed with
the center of a Delone polytope. We computed $241$ $T'$-generic
$T'$-secondary cones.  The new covering record is attained by a PQF
which has $1206$ orbits of Delone polytopes. Here we tried all
different centers. The new covering record was produced by using a 
center of a Delone polytope having the third largest symmetry group.

\subsection{New best known packing-coverings}

Finally we mention briefly a third application of our new theory.
Closely related to the lattice covering problem, is the
\textit{lattice packing-covering problem}.  It asks to minimize the
\textit{packing-covering constant}
\[
\gamma(L)=\frac{\mu(L)}{\lambda(L)}
\]
of a $d$-dimensional lattice $L$, where $\lambda(L)=\frac{1}{2}\min \{
\|\vec{v}\| : \vec{v} \in\Z^d\setminus\{\vec{0}\}\}$ is the so called
\textit{packing radius}.  For a detailed description of this problem
as well as for its interpretation as a convex optimization problem, we
refer the interested reader to \cite{sv-2005}.  So far, optimizing
over suitable $T$-secondary cones with respect to $4$-dimensional
subspaces $T$ of $\MS^7$, we found a new best known $7$-dimensional
packing-covering lattice with $\gamma(L)=1.499399\dots$.  The former
record holder was the lattice $\mathsf{E}^*_7$ with
$\gamma(\mathsf{E}^*_7) = \sqrt{7/3} = 1.527525\dots$. It
remains to undertake a systematic search for further best known, maybe
optimal packing-covering lattices.

%


\section*{Acknowledgements}

We like to thank Eva Bayer--Fluckiger and Gabriele Nebe for pointing
out application I, and we like to thank Yves Benoist for helpful
communications regarding Remark \ref{rem:closedgeodesics}. We thank
the anonymous referee for pointing out the connections of our paper to
previous works in algebraic geometry.
Regarding these we thank Gregory Sankaran for helpful comments.


\newcommand{\etalchar}[1]{$^{#1}$}
\providecommand{\bysame}{\leavevmode\hbox to3em{\hrulefill}\thinspace}

\end{document}